\documentclass[11pt]{article}
\usepackage[hmargin=0.9in,vmargin=1in]{geometry}
\usepackage[T1]{fontenc}
\usepackage[utf8]{inputenc}

\usepackage{lmodern,microtype,mathtools,amsthm,amssymb,xcolor,hyperref,bookmark}

\usepackage{thmtools}
\makeatletter
\@ifundefined{newcounteralias}{}{%
  \renewcommand\thmt@autorefsetup{%
    \@xa\def\csname\thmt@envname autorefname\@xa\endcsname
      \@xa{\thmt@thmname}%
  }%
}
\makeatother

\usepackage[shortlabels]{enumitem}
\usepackage{booktabs}
\usepackage{multicol}
\usepackage{xspace}
\usepackage{comment}
\hypersetup{
    unicode,
    colorlinks=true,
    citecolor=green!40!black,
    linkcolor=red!20!black,
    urlcolor=blue!40!black,
    filecolor=cyan!30!black
    }

\usepackage{cite}

\setlist[itemize]{leftmargin=2em,labelsep=.6em,topsep=\medskipamount,noitemsep}
\setlist[enumerate]{leftmargin=2em,labelsep=.5em,topsep=\medskipamount,noitemsep}

\newtheorem{theorem}{Theorem}[section]
\newtheorem{lemma}[theorem]{Lemma}
\newtheorem{proposition}[theorem]{Proposition}
\newtheorem{observation}[theorem]{Observation}
\newtheorem{corollary}[theorem]{Corollary}
\newtheorem{conjecture}[theorem]{Conjecture}

\newtheorem{claim}{Claim}[theorem]

\newcommand{\claimprefix}{}

\makeatletter
\newdimen\@savedtopsep
\renewenvironment{proof}[1][\proofname]{\par\@savedtopsep\topsep
\pushQED{\qed}\normalfont\topsep6\p@\@plus6\p@\trivlist
\item[\hskip\labelsep\itshape #1\@addpunct{.}]\topsep\@savedtopsep
\ignorespaces}{\popQED\endtrivlist\@endpefalse}
\makeatother

\newenvironment{claimproof}{\begin{proof}[Proof of Claim]}{\end{proof}}

\newcommand{\RR}{\mathbb{R}}

\newcommand{\cB}{\mathcal{B}}
\newcommand{\cC}{\mathcal{C}}
\newcommand{\cH}{\mathcal{H}}

\newcommand{\Oh}{\mathcal{O}}

\newcommand{\cG}{\mathcal{G}}
\newcommand{\cT}{\mathcal{T}}

\newcommand{\MWIS}{\textsc{MWIS}\xspace}
\newcommand{\MWISa}{\textsc{MWIS}$^*$\xspace}

\newcommand{\tw}{\textsf{tw}}
\newcommand{\treealpha}{\textsf{tree}\text{-}\alpha}
\newcommand{\yolov}{\textsf{tree}\text{-}\mu}

\newcommand{\wei}{\mathsf{w}}

\let\leq\leqslant\let\le\leq
\let\geq\geqslant\let\ge\geq
\renewcommand{\phi}{\varphi}
\newcommand{\rev}{\mathsf{rev}}

\usepackage{tikz}
\usetikzlibrary{patterns,decorations}

\usepackage[capitalise,nameinlink,sort]{cleveref}
\crefname{claim}{Claim}{Claims}
\Crefname{claim}{Claim}{Claims}
\crefname{enumi}{item}{items}
\crefname{problem}{Problem}{Problems}
\crefname{observation}{Observation}{Observations}
\crefname{conjecture}{Conjecture}{Conjectures}

\usepackage{authblk} 

\title{Induced matching treewidth in classes with a forbidden structure}
\author[1]{Valia Mitsou}
\author[2]{Paweł Rzążewski\thanks{Supported by the National Science Centre grant 2024/54/E/ST6/00094. A crucial part of the work was done when PRz was visiting Université Paris Cité}}
\author[1]{Giannos Stamoulis} 

\affil[1]{Université Paris Cité, CNRS, IRIF, UMR 8243, F-75013, Paris, France}
\affil[2]{Warsaw University of Technology}

\begin{document}
	\begin{titlepage}
	\date{}
	\maketitle
	
\begin{abstract}
    Induced matching treewidth is a tree-decomposition based parameter where we bound the size of an induced matching whose every edge intersects a single bag of the decomposition.
    It is known that many computational problems can be solved efficiently in graphs of bounded induced matching treewidth.

    In this paper we are interested in understanding what classes of graphs have bounded induced matching treewidth.
    We provide dichotomy theorems for all classes defined by excluding a single graph with respect to six basic containment relations: (induced) subgraph, (induced) topological minor, and (induced) minor.

    Furthermore, we study the case where the input graph has minimum degree at least 2, as well as the more restrictive case where it is 2-connected. 
    This can be assumed without loss of generality for many computational problems.
    It turns out that for classes excluding some graph as an induced minor, this restrictions yield some new cases of bounded induced matching treewidth.

    As a consequence, our results imply new algorithmic results, as well as the existence of balanced separators dominated by a constant number of vertices in the considered classes. This confirms some special cases of two conjectures of Gartland and Lokshtanov.
\end{abstract}

	\thispagestyle{empty}
\end{titlepage}

\tableofcontents
\thispagestyle{empty}
\newpage\setcounter{page}{1}

\section{Introduction}\label{sec:intro}
Graph decompositions are a general set of tools that allow us to break graphs into smaller pieces, connected to each other in somewhat restricted way.
One of the motivations behind this line of study is algorithmic: many computational problems can often be solved using some kind of dynamic programming over a decomposition,
and the ``complexity'' (usually referred to as \emph{width}) of the decompositions restricts the number of states one needs to consider.

Arguably, the best known type of graph decompositions are tree decompositions, discovered independently by several groups of authors with different motivationss~\cite{bertele1972nonserial,Halin1976,DBLP:journals/jct/RobertsonS84,Arnborg1987}.
Intuitively, a tree decomposition of a graph $G$ is a tree $\cT$ whose vertices (typically called \emph{nodes}) are equipped with subsets of the vertex set of $G$, called \emph{bags}, so that the general structure of $G$ is reflected by the structure of $\cT$. If each bag is ``simple,'' then the whole graph is considered ``simple.''
Perhaphs the most straightforward notion of simplicity is constant size.
This leads to the definition of \emph{treewidth}: a graph has treewidth $k$ if it admits a tree decomposition whose every bag is of size less than $k$.
Consequently, if a graph has small treewidth, then it ``resembles a tree.'' 

The drawback of treewidth is that its boundedness implies that the graph is sparse, i.e., it has a linear number of edges. Thus, this parameter is not suitable for studying classes of ``well-structured'' dense graphs, like chordal graphs.
As a remedy, in recent years there is a growing interest in \emph{structured tree decompositions} where we do not aim to obtain small bags, but rather bags of ``simple structure.''
A motivating example is the class of chordal graphs mentioned above: they are precisely graphs that admit tree decompositions whose every bag is a clique. Such cliques can be arbitrarily large, but their structure is indeed rather simple.

A generalization of this idea is the \emph{tree-independence number} (denoted by $\treealpha(G)$), introduced independently by Yolov~\cite{Yolov18} and by Dallard, Milani\v{c}, and \v{S}torgel~\cite{dallard2022firstpaper}. Here, for every bag of a tree decomposition, we mesure the size of a largest independent set in the graph induced by the bag. 
In particular, chordal graphs are precisely graphs of tree-independence number 1.
We remark that tree-independence number is always upper-bounded by treewidth (plus 1).

As already observed by Yolov~\cite{Yolov18}, tree-independence number might still be large on graphs of simple structure, like complete bipartite graphs. Thus, he introduced another parameter, \emph{minor-matching hypertreewidth}, which is now usually called \emph{induced matching treewidth} and denoted by $\yolov(G)$.
This time, for every bag of a tree decomposition, we measure the size of a largest induced matching whose every edge intersects the bag (note that the matching might not be contained in the bag).
It is not hard to observe that induced matching treewidth never exceeds tree-independence number, but can be much smaller. In particular, bicliques have induced matching treewidth 1.

When introducing a new graph parameter, there are three important questions to have in mind:
\begin{enumerate}[(Q1)]
    \item If the parameter is bounded, can we solve some interesting computational problems efficiently?\label{q1}
    \item If the parameter is bounded, can we compute a corresponding decomposition efficiently?\label{q2}
    \item Are there any interesting classes of graphs in which the parameter is bounded?\label{q3}
\end{enumerate}

Let us briefly discuss the state of art concerning these three questions for tree-independence number and induced matching treewidth.

Concerning~\ref{q1}, Yolov~\cite{Yolov18} showed that some classic problems including \textsc{Max Weight Independent Set} (\MWIS) and $k$-\textsc{Coloring} are solvable in polynomial time for graphs given with a tree decomposition witnessing that tree-independence number or induced matching treewidth is bounded.
Another family of problems that are tractable in the discussed setting concern packing induced subgraphs, possibly with distance constraints~\cite{dallard2022firstpaper,LMMORS26}.
Lima, Milani\v{c}, Mur\v{s}i\v{c}, Okrasa, Rzążewski, and \v{s}torgel~\cite{LMMORS26} showed that  $(\tw \leq r,\psi)$-\MWIS, a wide generalization of \MWIS (see \cref{sec:algos} for the definition) can be solved efficiently in classes of bounded tree-independence number (see also the refined version of this result~\cite{DBLP:journals/corr/abs-2601-15861}).
The authors conjectured that the same holds for induced matching treewidth, supporting this with several partial results.
The conjecture was recently confirmed by Bodlaender, Fomin, and Korhonen~\cite{BodlaenderFK26}.

Concerning~\ref{q2}, already in the first paper introducing both parameters, Yolov~\cite{Yolov18} showed that both tree-independence number and induced matching treewidth admit \textsf{XP}-approximation algorithm with respect to the value of the parameter. More precisely, given an $n$-vertex graph $G$ and an integer $k$, in time $n^{\Oh(k^3)}$ one can compute a tree decomposition witnessing that $\treealpha(G) = \Oh(k^3)$ (resp., $\yolov(G) = \Oh(k^3)$) or conclude that the corresponding parameter is larger than $k$.
Later, Dallard, Fomin, Golovach, Korhonen, and Milani\v{c}~\cite{dallard2022computing} proposed a refined algorithm that runs in time $2^{\Oh(k^2)} n^{\Oh(k)}$ and returns a tree decomposition showing that $\treealpha(G) \leq 8k$ or concludes $\treealpha(G) > k$. Using a structural observation of Yolov~\cite{Yolov18} one can use this algorithm to obtain analogous results for induced matching treewidth~\cite{LMMORS26}.

Interestingly, this approach is essentially optimal (up to improving the constant 8 and the $2^{\Oh(k^2)}$-factor). In particular, known reductions~\cite{dallard2022computing} exclude an exact algorithm even for $k=4$, as well as any \textsf{FPT}-approximation algorithm, under standard assumptions in complexity theory.

Classes considered in~\ref{q3} are typically defined by excluding some graphs under one of the standard containment relations: as a  subgraph, a topological minor, and a minor, as well as their induced counterparts, see~\cref{sec:prelim} for formal definitions of these relations.
For a relation $\star$ and a graph $H$, we say that a graph $G$ is \emph{$H$-$\star$-free} if it does not contain $H$ with respect to $\star$.

 Dallard, Milani\v{c}, and \v{S}torgel~\cite{dallard2022secondpaper} studied tree-independence number in classes defined by forbidding a single graph under the relations mentioned above. 
For each of the six considered relations, they fully characterized graphs $H$ whose exclusion bounds tree-independence number, see~\cref{tab:summary}.
Since then, the classes of small tree-independence number have been a topic of investigation~\cite{chudnovsky2025treewidth,abrishami2024tree,DBLP:journals/jctb/ChudnovskyHLS26,DBLP:conf/soda/ChudnovskyGHLS25}.

On the other hand, very little is known about which classes have bounded induced matching treewidth.
Let us mention that Abrishami, Briański, Czyżewska, McCarty, Rzążewski, and Walczak~\cite{AbrishamiBCMMRW25} showed that induced matching treewidth and tree-independence number are functionally equivalent for graphs that exclude some fixed biclique as an induced subgraph. However, for classes that contain arbitrarily large bicliques, we are not aware of any results.

\subsection*{Our results.}
In this paper we initite a systematic study of classes of graphs that have bounded induced matching treewidth,
see~\cref{tab:summary}. The definitions of the graphs listed throughout this section can be found in~\cref{sec:prelim}.

\begin{table}[t]
    \centering    
    {\small
\begin{tabular}{l||c | c |c}
    \toprule 
    Relation & $\treealpha$~\cite{dallard2022secondpaper} & $\yolov$ & $\yolov$ when $\delta \geq 2$ \\ \midrule \midrule
     Subgraph &\multicolumn{3}{c}{  disjoint union of subdivided claws   } \\ \midrule
     Induced subgraph & $P_3$ or edgeless    &  \multicolumn{2}{c}{$P_3$ or a matching} \\ \midrule
     Topological minor &  \multicolumn{3}{c}{  planar and subcubic   }  \\ \midrule
     Minor &   \multicolumn{3}{c}{ planar    }  \\           \midrule
     Induced topological minor &  $C_4$, $K_4^-$, or edgeless   &  \multicolumn{2}{c}{$C_4$, $K_4^-$, 1-tipi, or a matching} \\ \midrule
     Induced minor & $K_{2,q}$ for some $q$, $W_4$,   &  $K_{2,q}$ for some $q$, $W_4$, & $K_{2,q}$ for some $q$, $W_4$, \\ 
                    & $K_5^-$, or edgeless    &  $K_5^-$, a tipi, or a matching &    $K_5^-$, or a windmill \\
     \bottomrule
\end{tabular}
}
    \caption{Maximal (w.r.t. a given containment relation) graphs whose exclusion (w.r.t. that relation) implies boundedness of tree-independence number or induced matching treewidth. In the last column we additionally assume that graphs have no vertices of degree at most 1.}
    \label{tab:summary}
    \end{table}

\paragraph{Dichotomy theorems.}    
First, we show results analogous to ones obtained for tree-independence number by Dallard, Milani\v{c}, and \v{S}torgel~\cite{dallard2022secondpaper}:
for each of the six basic containment relations, we fully characterize which classes defined by excluding a single graph have bounded induced matching treewidth.

For subgraphs, induced subgraphs, minors, and topological minors, the dichotomies almost coincide with the ones for tree-independence number. Note that graphs with no large induced matching have trivially bounded induced matching treewidth, just like graphs with no large independent set have trivially bounded tree-independence number. 

\begin{theorem}\label{thm:dichoeasy}
    Let $H$ be a graph with at least 2 vertices.
    The following hold.
    \begin{enumerate}[(1.)]
        \item The class of $H$-subgraph-free graphs has bounded induced matching treewidth if and only if $H$ every component of $H$ is a subcubic tree with at most one vertex of degree 3. \label{it:easy-subgraph}
        \item The class of $H$-free graphs has bounded induced matching treewidth if and only if $H$ is $P_3$ or a graph with maximum degree 1.\label{it:easy-inducedsubgraph}
        \item The class of $H$-minor-free graphs has bounded induced matching treewidth if and only if $H$ is planar.\label{it:easy-minor}
        \item The class of $H$-topological-minor-free graphs has bounded induced matching treewidth if and only if $H$ is subcubic and planar.
        \label{it:easy-topminor}
    \end{enumerate}
    All unbounded cases remain so, even when restricted to 2-connected graphs.
\end{theorem}

\begin{figure}[t]
    \centering    
    \includegraphics[page=1]{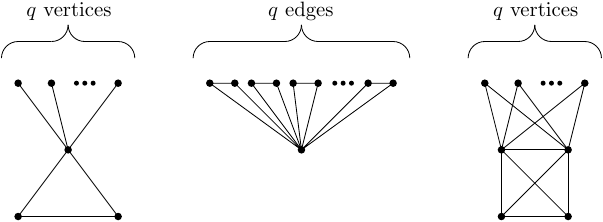}
    \caption{The $q$-tipi (left), $q$-windmill (middle), and $\widehat{K}_4^{(q)}$ (right).}
    \label{fig:tipiwindmillfullhouse}
\end{figure}

The landscape gets a bit richer when we exclude a graph as an induced topological minor: here we get a new bounded case: \emph{1-tipi}, also known as the \emph{paw}, see~\cref{fig:tipiwindmillfullhouse}.

\begin{restatable}{theorem}{thmdichosubdivision}
\label{thm:dichosubdivision}
    Let $H$ be a graph with at least 2 vertices.
    The class of $H$-topological-minor-free graphs has bounded induced matching treewidth if and only if $H$ is not contained in one of the following graphs: $C_4$, $K_4^-$, 1-tipi, or a matching as an induced topological minor.
    All unbounded cases remain so, even when restricted to 2-connected graphs.
\end{restatable}

For $H$-induced-minor-free graphs, we actually get an infinite family of bounded cases: it turns out that excluding \emph{any tipi} (i.e., the $q$-tipi for any $q$) bounds the induced matching treewidth.

\begin{restatable}{theorem}{thmdichoindminor}
\label{thm:dichoindminor}
    Let $H$ be a graph with at least 2 vertices.
    The class of $H$-induced-minor-free graphs has bounded induced matching treewidth if and only if $H$ is an induced minor of one of the following graphs: $K_{2,q}$ for some $q\in\mathbb{N}$, $W_4$, $K_5^-$, a tipi, or a matching.
\end{restatable}

\paragraph{Restricting to graphs of minimum degree 2.} 
Recall that the motivation to consider induced matching treewidth was that tree-independence number is unbounded in some very simple classes of graphs, like bicliques.
Admittedly, induced matching treewidth has similar drawback. Consider for example the graph obtained from the biclique with each part of size $n$ by adding a private neighbor to every vertex.
As noted by Yolov~\cite{Yolov18} this graph has induced matching treewidth $n$.

However, intead of defining a new parameter, we can take another route.
If our motivation is algorithmic, perhaphs degree-1 vertices do not play a crucial role and can be removed by some preprocessing rule.
In general, we can try to apply various (problem-dependent) preprocessing rules to the considered graph and only then ask if the considered parameter is bounded.
This idea is not new, but previously it was applied to other (not tree-decomposition-based) graph parameters, like \emph{clique-width}~\cite{GaspersHP19,DabrowskiMNPR23} and \emph{mim-width}~\cite{DBLP:journals/corr/abs-2004-05022}.

We investigate this direction for tree-independence number and induced matching treewidth in classes excluding a fixed subtructure.
First, let us restrict our attention to graphs with minimum degree 2 as for most classic computational problems we can indeed easily handle vertices of degree at most 1.
It turns out that for tree-independence number this change in the setting does not alter the dichotomies for any of the six consider containment relations.
Similarly, for induced matching treewidth and (induced) subgraph, (induced) topological minors, and minors,
the classification remains the same as before.

However, for $H$-induced-minor-free graphs, we get new family of bounded cases: if $H$ is a \emph{windmill}, i.e., the $q$-windmill for any $q$; see~\cref{fig:tipiwindmillfullhouse}.

\begin{restatable}{theorem}{thmdichoindminordeg}
\label{thm:dichoindminordeg}
    Let $H$ be a graph.
    The class of $H$-induced-minor-free graphs, restricted to graphs with minimum degree at least 2,
    has bounded induced matching treewidth if and only if $H$ is an induced minor of one of the following graphs: $K_{2,q}$ for some $q\in\mathbb{N}$, $W_4$, $K_5^-$, or a windmill.
\end{restatable}

Graphs that exclude a fixed windmill as an induced minor were considered by Bonnet, Duron, Geniet, Thomass\'e, and Wesolek~\cite{BonnetDGTW26}; we will get back to this soon.

\paragraph{Restricting to 2-connected graphs.} 
We can further generalize the idea of removing vertices of degree 1:
we can assume that the input graph is 2-connected.
Indeed, for many problems such an assumption is justified (see e.g.~\cref{sec:algos}).
Note that this setting generalizes the one considered before, as every 2-connected graph has no vertices of degree at most 1.
Again, it turns out that for two parameters and six containment relations studied in this paper,
only the case of induced matching treewidth of $H$-induced-minor-free graphs differs from the basic dichotomies without any further restrictions imposed on the graph.

In this setting, we get one more bounded case: $\widehat{K}^{(1)}_4$, for simplicity denoted as $\widehat{K}_4$, see~\cref{fig:tipiwindmillfullhouse}.

\begin{restatable}{theorem}{thmindminortwoconn}
\label{thm:indminortwoconn}
    Let $H$ be a graph.
    The class of $H$-induced-minor-free graphs, restricted to 2-connected graphs,
    has bounded induced matching treewidth if $H$ is an induced minor of one of the following graphs: $K_{2,q}$ for some $q$, $W_4$, $K_5^-$, a windmill, a matching, or $\widehat{K}_4$.
    If $H$ is not one of these graphs, nor  $\widehat{K}_4^{(q)}$ for some fixed $q$, then 2-connected $H$-induced-minor-free graphs have unbounded induced matching treewidth.
\end{restatable}

Our \cref{thm:indminortwoconn} leaves one family of open cases: $H=\widehat{K}_4^{(q)}$ for some $q \geq 2$.
For them, we were not able to prove either boundedness nor unboundedness of induced matching treewidth for 2-connected $H$-induced-minor-gree graphs. We believe in the former outcome and we state the following conjecture.

\begin{restatable}{conjecture}{conjfullhouses}
    \label{conj:fullhouses}
    For every fixed $q \geq 1$, the class of 2-connected $\widehat{K}_4^{(q)}$-induced-minor-free graphs has bounded induced matching treewidth.
\end{restatable}

\paragraph{Algorithmic consequences.}
    Now, let us turn our attention to algorithmic implications of our results.
    We consider two general families of problems:  $(\tw \leq r,\psi)$-\MWIS and \textsc{Max Induced $H$-Colorable Subgraph}.
    The first one asks for a maximum-weight induced subgraph of treewidth at most $r$, satisfying given \textsf{CMSO}$_2$ formula,
    while the second one asks for a maximum-weight induced subgraph homomorphic to a fixed simple graph $H$; see~\cref{sec:algos} for precise definitions.
    Let us mention that the first family captures, e.g., \MWIS, \textsc{Feedback Vertex Set}, and \textsc{Even Cycle Transversal},
    while the second family captures, e.g., \MWIS (again), \textsc{Odd Cycle Transversal}, and \textsc{List $k$-Coloring} for any fixed $k$.

    We show that for all problems described by the formalisms above, it is indeed enough to solve instances that are 2-connected. Thus, combining this with known algorithms for graphs of bounded induced matching treewidth~\cite{Yolov18,BodlaenderFK26} and our \cref{thm:indminortwoconn}, we obtain the following algorithmic corollary.
    (The first column follows from bounded tree-independence number and thus, it was known before, while the second column contains new results.)

        \begin{corollary}\label[corollary]{cor:algsummary}
        Every problem that can be described as \textsc{Max Induced $H$-Colorable Subgraph} for some fixed simple graph $H$ or $(\tw \leq r,\psi)$-\MWIS for fixed $r$ and a formula $\psi$, can be solved in polynomial time when restricted to graphs from the following classes:
    \begin{multicols}{2}
    \begin{itemize}
        \item $K_{2,q}$-induced-minor-free for any fixed $q$,
        \item $W_4$-induced-minor-free,
        \item $K_5^-$-induced-minor-free,
        \item $q$-windmill-induced-minor-free for any fixed $q$,
        \item $\widehat{K}_4$-induced-minor-free.
    \end{itemize}
    \end{multicols}
\end{corollary}

Let us mention some context of these results. 
First, as already mentioned,  Bonnet, Duron, Geniet, Thomass\'e, and Wesolek~\cite{BonnetDGTW26} considered graphs that exclude a fixed windmill as an induced minor, and showed that \MWIS can be solved in polynomial time in every such class.
Our~\cref{cor:algsummary} not only heavily generalizes this result, but also shows the underlying reason \emph{why} this class admits algorithms for classic problems.

Second, our results confirm special cases of the following conjecture of Gartland and Lokshtanov~\cite{gartland2023quasi}.

\begin{conjecture}[Gartland and Lokshtanov]\label[conjecture]{conj:GLalg}
    For any fixed planar graph $H$, integer $r$, and a formula $\psi$,
    $(\tw \leq r,\psi)$-\MWIS and \textsc{List $3$-Coloring} can be solved in polynomial time in $H$-induced-minor-free graphs.
\end{conjecture}

Very recently, \cref{conj:GLalg} was disproved by Bonnet and Chang\cite{BonnetChang26} for $H$ being the $5 \times 5$ grid.
However, the statement was confirmed for some families of simple planar graphs $H$~\cite{DBLP:conf/soda/LokshantovVV14,DBLP:journals/talg/GrzesikKPP22,DBLP:journals/siamcomp/AbrishamiCPRS24,DBLP:journals/talg/ChudnovskyMPPR26,DBLP:journals/algorithmica/ChudnovskyHSZ21,DBLP:conf/esa/GalbyLMN25,DBLP:journals/algorithmica/HoangKLSS10,DBLP:journals/algorithmica/KlimosovaMMNPS20,DBLP:journals/combinatorica/BonomoCMSSZ18} and it remains a very interesting open problem to understand for which graphs $H$ the conjecture actually holds.
Even its relaxed variants, asking for quasipolynomial-time algorithms, or at least a quasipolynomial-time approximation scheme (QPTAS) for \MWIS, remain very challenging and are shown only in restricted cases~\cite{DBLP:conf/focs/GartlandL20,gartland2021finding,DBLP:journals/jctb/NguyenSS24,DBLP:journals/jctb/BonamyBDEGHTW24,DBLP:conf/swat/BonnetCMPR26,thetas,chudnovsky2025dominated,BonnetDGTW26,DBLP:conf/sosa/PilipczukPR21}.

\paragraph{Dominated balanced separators.}
Finally, we show one more consequence of our results.
For a constant $t$, a class of graphs is said to admit \emph{$t$-dominated balanced sepatator} 
if every graph from the class a set ot at most $t$ vertices,
such that removing these vertices along with their neighbors breaks the graph into small pieces; see~\cref{sec:dbs} for precise definitions.
We show that boundedness of induced matching treewidth, even for just 2-connected members of the considered class, yields $\Oh(1)$-dominated balanced separators. Thus we obtain the following.

    \begin{restatable}{corollary}{cordbs}\label[corollary]{cor:dbs}
    The following classes of graphs admit $t$-dominated balanced separators, where $t$ is a constant depending only on the class.
    \begin{multicols}{2}
    \begin{itemize}
        \item $K_{2,q}$-induced-minor-free for any fixed $q$,
        \item $W_4$-induced-minor-free,
        \item $K_5^-$-induced-minor-free,
        \item $q$-windmill-induced-minor-free for any fixed $q$,
        \item $\widehat{K}_4$-induced-minor-free.
    \end{itemize}
    \end{multicols}
\end{restatable}

\cref{cor:dbs} confirms special cases of another conjecture by Gartland and Lokshtanov~\cite{gartland2023quasi}.
\begin{conjecture}[Gartland and Lokshtanov]\label[conjecture]{conj:GLdbs}
    For any fixed planar graph $H$ there exists $t$ such that $H$-induced-minor-free graphs admit $t$-dominated balanced separators.
\end{conjecture}

\cref{conj:GLdbs} is also very challenging and, despite interest in the community, was only confirmed for very restricted graphs $H$~\cite{BacsoLMPTL19,chudnovsky2024quasi,chudnovsky2025dominated,thetas}.

\paragraph{Organization of the paper.}
In \cref{sec:prelim} we define main notions used in the paper and show some basic results.
In \cref{sec:unbounded} we construct some families of graphs with unbounded induced matching treewidth.
In \cref{sec:tipi} we discuss graphs that exclude a tipi as an induced minor.
In \cref{sec:dichos} we show \cref{thm:dichoeasy,thm:dichosubdivision,thm:dichoindminor}, i.e., dichotomy theorems for classes that exclude one fixed graph.
In \cref{sec:noleaves} we dicuss the problem where we additionally assume that th graphs considered have minimum degree 2, and prove~\cref{thm:dichoindminordeg}.
In \cref{sec:2connected} we restrict ourselves to 2-connected graphs wne show~\cref{thm:indminortwoconn}.
In \cref{sec:algos} we discuss algorithmic consequences of our results and show~\cref{cor:algsummary}.
In \cref{sec:dbs} we discuss dominated balanced separators an show~\cref{cor:dbs}.
Finally, in \cref{sec:outro}, we point out some directions for further research.

\newpage
\section{Definitions and basic tools}\label{sec:prelim}
For an integer $n$, by $[n]$ we denote the set $\{1,2,\ldots,n\}$.
For a set $X$, by $2^X$ we denote the family of all subsets of $X$.

\paragraph{Graph theory.}
All graphs in the paper are undirected and simple, i.e., without loops or parallel edges.
Let $G$ be a graph and consider $X \subseteq V(G)$.
By $G[X]$ we denote the subgraph of $G$ \emph{induced} by $X$ and by $G-X$ we denote $G[V(G) \setminus X]$.
We use $N[X]$ for the closed neighborhood of $X$, i.e., the set consisting of $X$ and all vertices adjacent to a vertex of $X$.
If it does not lead to confusion, sometimes we identify induced subgraphs with their vertex sets.

We say that a graph is \emph{biconnected} if it is 2-connected or isomorphic to $K_2$.
A \emph{block} is an inclusion-wise biconnected (induced) subgraph of a graph.
A \emph{block-cut tree} of a connected graph $G$ is the bipartite graph whose vertices are blocks and cutvertices of $G$,
and a cutvertex $v$ and a block $b$ are adjacent if and only if $v \in b$.
It is well-known that the block-cut tree is indeed a tree. Note that all leaves of this tree are blocks.
We treat block-cut trees as rooted trees, where the root $r$ is an arbitrarily chosen leaf and which defines a parent-child relation between adjacent nodes.
For a block node $b \neq r$, if $v_b$ is the cutvertex that is the parent of $b$ in $T$, then the parent of $v_b$ in $T$ is called the \emph{parent block}.
Analogously we define \emph{children blocks} and \emph{leaf blocks}.

A \emph{matching} in a graph is a set of pairwise disjoint edges.
For a matching $M$, by $V(M)$ we denote the set of vertices of edges in $M$.
A matching $M$ in a graph $G$ is \emph{induced} if the edge set of $G[V(M)]$ is $M$, i.e., there are no edges joining vertices belonging to two distinct edges in $M$.
Given a graph $G$ and two vertex subsets $A,B\subseteq V(G)$, an \emph{induced $(A,B)$-matching} of $G$ is an induced matching $M$ of $G$ whose edges have one endpoint in $A$ and one endpoint in $B$.

\paragraph{Graph containment relations.}
For a graph $G$ and an edge $xy$, \emph{contracting} $xy$ is the operation of introducing a new vertex adjacent to every vertex adjacent to $x$ or $y$, and then removing $x$ and $y$. 
A contraction is \emph{topological} if one of the vertices $x,y$ is of degree 2 in $G$. For a graph $G$ and an edge $xy$, \emph{subdividing} $xy$ is the operation of adding a new vertex adjacent to $x$ and $y$.

For a containment relation $\star$, we say that a graph $G$ contains a graph $H$ as $\star$, if $H$ can be obtained from $G$ by a sequence of allowed \emph{elementary operations}, according to the following table:

\medskip
\begin{tabular}{l|| c | c | c | c}
    \toprule 
             & vertex & edge & edge & topological \\
    Relation & deletion & deletion & contraction & contraction \\ \midrule \midrule
     Subgraph & + & +  & - & - \\
     Induced subgraph & + & -  & - & - \\
     Topological minor & + & +  & - & + \\
     Minor & + & +  & + & + \\
     Induced topological minor & + & -  & - & + \\
     Induced minor & + & -  & + & + \\
     \bottomrule
\end{tabular}
\medskip

We will interchangeably use terms \emph{induced topological minor} and \emph{induced subdivision}. 

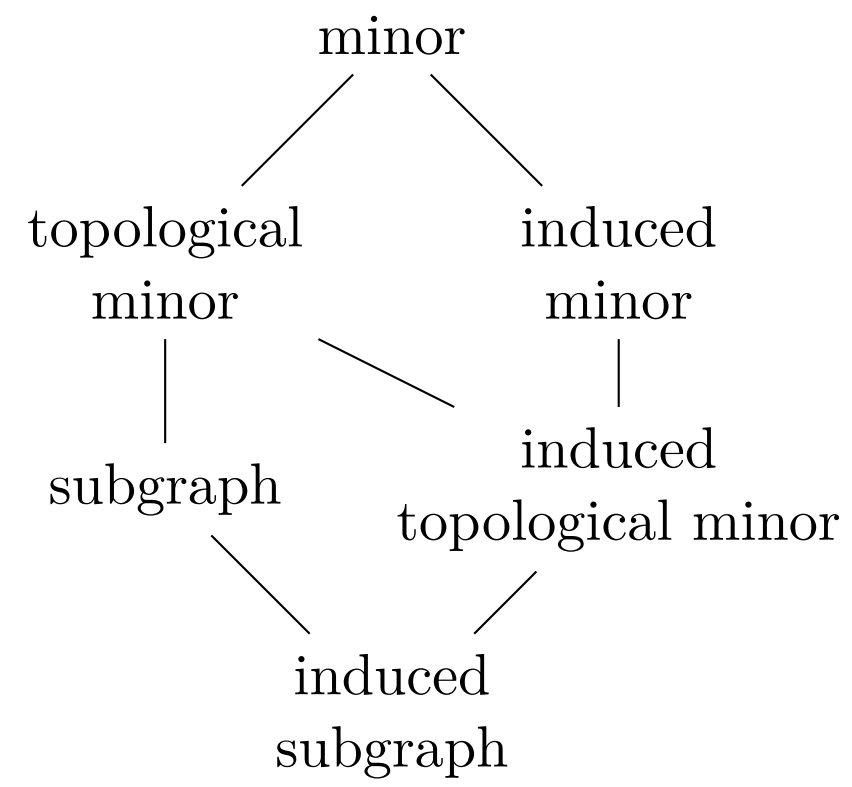
\begin{figure}[t]
    \centering    
    \begin{tikzpicture}[every node/.style={draw = none, fill = none},scale=0.5]

        \node[align=center] (minor) at (0,0) {minor};
        \node[align=center] (topminor) at (-3,-3) {topological\\minor};
        \node[align=center] (indminor) at (3,-3) {induced\\minor};
        \node[align=center] (subgraph) at (-3,-6) {subgraph};
        \node[align=center] (indsubdiv) at (3,-6) {induced\\ topological minor};
        \node[align=center] (indsubgraph) at (0,-9) {induced\\subgraph};

        \draw (minor) -- (topminor);
        \draw (minor) -- (indminor);

        \draw (topminor) -- (subgraph);
        \draw (topminor) -- (indsubdiv);
        \draw (indminor) -- (indsubdiv);

        \draw (subgraph) -- (indsubgraph);
        \draw (indsubdiv) -- (indsubgraph);
    \end{tikzpicture}
    \caption{The considered containment relations as a poset.
    A relation $\star$ is greater than (above) a relation $\ast$ if every graph $G$ that contains $H$ as $\ast$ also contains $H$ as $\star$. In other words, $H$-$\star$-free graphs  form a subclass of $H$-$\ast$-free graphs.    
    }    
    \label{fig:relations}
\end{figure}

Note that all relations discussed above are transitive, i.e., if $G$ contains $H$ and $H$ contains $F$, then $G$ contains $F$ (with respect to the considered relation).
Furthermore, the relations between these relations are depicted in~\cref{fig:relations}

A class is \emph{hereditary} if it is closed under vertex deletion (i.e., taking induced subgraphs).
Note that if $\star$ is any of the relations discussed above and $H$ is any graph,
then $H$-$\star$-free graphs form a hereditary class.

\paragraph{Some special graphs.}

For graphs $H_1,H_2$, by $H_1+H_2$ we denote the graph obtained as a disjoint union of $H_1$ and $H_2$.

By $P_k$ (resp., $C_k$) we denote the path (resp., the cycle) with $k$ vertices.
We say that a graph is a \emph{linear forest} if it is isomorphic to a disjoint union of paths.
If it does not lead to confusion, a linear forest with every component having at most two vertices will be called a \emph{matching}.

A \emph{wheel} is a graph obtained from a cycle by adding a universal vertex.
The \emph{4-wheel}, denoted by $W_4$, is the wheel obtained from $C_4$.

For $k \in\mathbb{N}$, the \emph{$k$-tipi} is the graph that is obtained from $K_3$ after attaching $k$ vertices of degree one to a single vertex of $K_3$; see~\Cref{fig:tipiwindmillfullhouse}.
For $k\in\mathbb{N}$, the \emph{$k$-windmill} is the graph that is obtained from a matching with $k$ edges by adding a universal vertex, see~\cref{fig:tipiwindmillfullhouse}.
Note that the $k$-tipi is an induced subgraph of the $(k+1)$-windmill.
We say that $G$ is a \emph{tipi} (resp., a \emph{windmill}) if it is a $k$-tipi (resp. $k$-windmill) for some $k\in\mathbb{N}$.

The \emph{bull} is the graph obtained from $K_3$ by adding two degree-1 vertices, each adjacent to a different vertex of the triangle.
The \emph{house} is the graph obtained from $C_4$ by adding a new vertex, adjacent to two consecutive vertices of the cycle.

By $K_k$ we denote the complete graph on $k$ vertices.
By $K_k^-$ we denote the graph obtained from $K_k$ by deleting a single edge.
For integers $k_1,k_2,\ldots,k_\ell$, by $K_{k_1,k_2,\ldots,k_\ell}$ we denote the complete $\ell$-partite graph with parts of size, respectively, $k_1,k_2,\ldots,k_\ell$.

For $q \geq 1$, by $\widehat{K}_4^{(q)}$ we denote the graph obtained from $K_4$ by adding $q$ new, pairwise nonadjacent vertices, all adjacent with the same pair of vertices of the $K_4$. We shortly write $\widehat{K}_4 \coloneqq \widehat{K}_4^{(1)}$.

A \emph{$k$-wall} is the $k \times k$ hexagonal grid.
A \emph{subdivided $k$-wall} is any graph obtained from the $k$-wall by repeatedly subdividing edges.
A \emph{properly subdivided $k$-wall} is obtained from the $k$-wall by subdividing each edge at least once.
Note that every (subdivided) wall is subcubic, and so is the line graph of each properly subdivided wall. 

For graphs $G$ and $H$, let $G \odot H$ denote the graph obtained from $G$ by adding $V(G)$ copies of $H$, one for each vertex $v$ of $G$,
and adding all edges between $v$ and the vertices of its copy of $H$.

\paragraph{Treewidth, tree-independence number, and induced matching treewidth.}
Let $G$ be a graph.
A \emph{tree decomposition} of $G$ is a pair $\cT = (T, \beta)$ where $T$ is a tree and $\beta : V(T) \to 2^{V(G)}$ is a \emph{bag} function that satisfies the following two properties:
\begin{itemize}
    \item for each vertex $v \in V(G)$, the set $\{t \in V(T) ~|~ v \in \beta(t) \}$ is nonempty and induces a subtree of $T$,
    \item for each edge $e \in E(G)$ there is $t \in V(T)$ such that $e \subseteq \beta(t)$.
\end{itemize}
A tree decomposition $(T,\beta)$, where $T$ is a path, is called a \emph{path decomposition}.

Let $\cT = (T, \beta)$ be a tree decomposition of a graph $G$.
The \emph{width} of $\cT$ is defined as $\max_{t \in V(T)} |\beta(t)|-1$.
The \emph{independence number} of $\cT$ is the size of a largest independent set in $G[\beta(t)]$, over all $t \in V(T)$
The \emph{induced matching width} of $\cT$ is the largest size of an induced $(\beta(t), V(G))$-matching,
over all $t \in V(T)$

The \emph{treewidth} (resp.  \emph{tree-independence number},  \emph{induced matching treewidth})
of $G$, denoted by $\tw(G)$ (resp., $\treealpha(G)$, $\yolov(G)$) is the minimum width (resp., independence number, induced matching width) of a tree decomposition of $G$.

The following proposition summarizes some known relations between these three parameters, see~\cite{Yolov18,AbrishamiBCMMRW25,LMMORS26}.

\begin{proposition}\label[proposition]{prop:properties}
    The following hold for every graph $G$.
    \begin{enumerate}[(1.)]
        \item $\yolov(G) \leq \treealpha(G) \leq \tw(G)+1$, \label{it:yolov-tin-tw}
        \item $\treealpha(G) = \yolov(G\odot K_1)$, \label{it:yolov-odot}
        \item $\tw(G) \leq 2\yolov(G) \cdot \Delta(G)^2(\Delta(G)+1)$, where $\Delta(G)$ is the maximum vertex degree in $G$. \label{it:yolov-Delta}
    \end{enumerate}
\end{proposition}

Let us also show one simple property of tree-independence number.

\begin{lemma}\label[lemma]{lem:block-reduction}
    Let $\cG$ be a hereditary class, and let $\cG_{\mathrm{2\textnormal{-}conn}}$ consist of the 2-connected graphs in $\cG$.
    Then $\cG$ has bounded tree-independence number if and only if $\cG_{\mathrm{2\textnormal{-}conn}}$ has bounded tree-independence number.
\end{lemma}
\begin{proof}
    The forward implication is immediate. Conversely, suppose that every graph in $\cG_{\mathrm{2\textnormal{-}conn}}$ has tree-independence number at most $k \geq 1$.

    Pick any $G \in \cG$. We may assume that $G$ is connected, as otherwise we can easily obtain a tree decomposition of $G$ by combining tree decompositions of its components.
    Every block of $G$ is either 2-connected or isomorphic to $K_2$;
    hence, since $\cG$ is hereditary, it has a tree decomposition of independence number at most $k$.

    Fix such a tree decomposition $\cT_B$ of each block $b$ of $G$.
    For each cutvertex $v$ of $G$, we create a node $t_v$ whose bag is $\{v\}$.
    Now we combine these decompositions along the block-cut tree of $G$.
    More precisely, if a block $b$ contains a cutvertex $v$, we join the node $t_v$ to an arbitrary node of $\cT_b$ whose bag contains $v$.
    It is straightforward to verify that this way we obtain a tree decomposition of $G$ of independence number at most $k$.
\end{proof}

\paragraph{BFS-layerings.}
For a connected graph $G$ and a vertex $r$, the \emph{BFS-layering of $G$ rooted at $r$} is the partition $( L_0,\ldots,L_q)$ of $V(G)$, where $L_0=\{r\}$ and for $i \geq 1$, the layer $L_i$ consists of all neighbors of vertices in $L_{i-1}$ that do not belong to $\bigcup_{j\le i}L_{j}$.
For a vertex $v$ in layer $L_i$, for some $i \in [0,q]$, a \emph{root path} is an $r$-$v$-path visiting each layer $0,1,\ldots,i$ exactly once; clearly such a path exists for any vertex $v$.

The following lemma shows how to turn BFS-layerings into tree decompositions of small induced matching width.

\begin{lemma}\label[lemma]{lem:from-bfs-to-tree}
    Let $G$ be a connected graph, and consider its BFS layering $(L_0,L_1,\ldots,L_q)$ rooted at some vertex.
    If, for each $i \geq 1$, it holds that:
    \begin{itemize}
        \item a largest induced matching in $G[L_i]$ has at most $a$ edges and
        \item a largest induced $(L_{i-1},L_i)$-matching has at most $b$ edges,
    \end{itemize}
    then $\yolov(G) \leq 2a+3b$.
\end{lemma}    
\begin{proof}
    Consider a path decomposition with bags of consecutive vertices 
    \[L_0 \cup L_1, L_1 \cup L_2, \ldots, L_{q-1} \cup L_q.\]
    It is straightforward to verify that it is indeed a path decomposition.
    Fix any bag $X= L_{i-1} \cup L_i$ and let $M$ be an induced matching whose every edge intersects $X$.
    We can partition the edges of $M$ into five groups: edges between $L_{i-2}$ and $L_{i-1}$ (if $L_{i-2}$ exists), edges in $L_{i-1}$, edges between $L_{i-1}$ and $L_{i}$, edges in $L_i$, and edges between $L_{i}$ and $L_{i+1}$ (if $L_{i+1}$ exists).
    The bound follows from upper-bounding the number of edges of each type independently.
\end{proof}

\paragraph{Induced minors.}

Let $G$ and $H$ be graphs. 
An \emph{induced minor model} of $H$ in $G$ is a function $\phi: V(H)\to 2^{V(G)}$ such that for every $v\in V(H)$ the subgraph $G[\varphi(v)]$ is connected and also $uv\in E(H)$ if and only if there is at least one edge in $G$ between $\phi(u)$ and $\phi(v)$. The sets $\phi(v)$ for $v\in V(H)$ are called \emph{branch sets}.

A walk $W=v_1,\ldots,v_\ell$ in a graph $G$ is a \emph{thin walk} if all its vertices are pairwise distinct, except possibly from $v_1 \ v_\ell$,
and all vertices are of degree 2 in $H$.  Vertices $v_2,\ldots,v_{\ell-1}$ are called \emph{internal vertices} of the walk. 

The following lemma is a slight generalization of the observation of Milani\v{c}, Piva\v{c}~\cite[3.2]{DBLP:journals/corr/abs-2405-15543} that is given by the same proof.

\begin{lemma}\label[lemma]{lem:minimal}
    Let $H$ be a graph, $W$ be a set of edge-disjoint thin walks in $H$, and $U$ be the set consisting of all internal vertices of walks in $W$ and all vertices of degree at most 1 in $H$.
    Let $H$ be an induced minor of a graph $G$.
    Then, there exists an induced minor model $\phi : V(H) \to 2^{V(G)}$ of $H$ in $G$ such that $|\phi(u)|$ = 1 for all $u \in U$.
\end{lemma}

\begin{corollary}\label[corollary]{cor:linearforests}
    If $H$ is a linear forest, then the class of $H$-free graphs coincides with the class of $H$-induced-minor-free graphs.
\end{corollary}

\newpage
\section{Families of unbounded induced matching treewidth}\label{sec:unbounded}

We start this section by listing a few families of graphs with unbounded induced matching treewidth; the definitions of the considered graphs were introduced in \cref{sec:prelim}.

\begin{lemma}\label[lemma]{lem:constructions}
    The following graphs have induced matching treewidth at least $k$:
    \begin{enumerate}[(1.)]
        \item any subdivided $72k$-wall, \label{it:wall}
        \item the line graph of any properly subdivided $72k$-wall, \label{it:linewall}
        \item the graph $M_k$ constructed from two $k$-edge matchings, by adding all edges between one matching and another, \label{it:Mk}
        \item the graph $K_{k,k}^\circ$ obtained from $K_{k,k}$ by fixing a perfect matching and, for each matching edge $xy$, adding a new 5-vertex path with endpoints $x$ and $y$,  \label{it:Kkkcirc}
        \item $K_{k,k} \odot K_2$, \label{it:KkkodotK2} 
        \item $K_{k,k} \odot K_1$. \label{it:KkkodotK1}               
    \end{enumerate}    
    Furthermore, for $k \geq 2$, the graphs constructed in items \ref{it:wall}--\ref{it:Kkkcirc}  of \cref{lem:constructions} are 2-connected, while graphs constructed in item \ref{it:KkkodotK2} have minimum degree 2.
\end{lemma}
\begin{proof}
    The first two items follow from recalling that (subdivided) $k$-walls and their line graphs have treewidth $k$.
    Since they are subcubic graphs, the bound on induced matching treewidth follows immediately from \cref{prop:properties} \ref{it:yolov-Delta}.

    For the third item, recall a well-known property of tree decompositions that for every tree decomposition $\cT$ there exists a vertex whose whole closed neighborhood is contained in a single bag.
    The neighborhood of every vertex in $M_k$ contains a $k$-edge induced matching, proving that $\yolov(M_k) \geq k$.

    For the last three items, note that all these graphs contain $K_{k,k} \odot K_1$ as an induced subgraph.
    Thus, the claim follows from  \cref{prop:properties}~\ref{it:yolov-odot} and the facts that $\treealpha(K_{k,k}) \geq k$~\cite{dallard2022firstpaper} and that induced matching treewidth does not increase when deleting vertices.    
\end{proof}

Let us show one more property of induced minors of $K_{k,k} \odot K_2$.

\begin{lemma}\label[lemma]{lem:2connbiclique}
    Let $H$ be a 2-connected graph with at least 4 vertices.
    Then $H$ is an induced minor of $K_{k,k}$ if an only if it is an induced minor of $K_{k,k} \odot K_2$.   
\end{lemma}
\begin{proof}
    \renewcommand{\claimprefix}{\ref{lem:2connbiclique}}
    \setcounter{claim}{0}

    The forward implication is obvious, as $K_{k,k}$ is an induced subgraph of $K_{k,k} \odot K_2$. 
    Now, suppose that $K_{k,k} \odot K_2$ contains $H$ as an induced minor and consider an induced minor model $\phi : V(H) \to 2^{V(K_{k,k} \odot K_2)}$.

    Denote the vertex set of the biclique $K_{k,k}$ in $K_{k,k} \odot K_2$ by $B$ and let $Z = V(K_{k,k} \odot K_2) \setminus B$.
    
    \begin{claim}\label{clm:onlybiclique}
        Let $z \in Z$ be such that $z \in \phi(v)$ for some $v \in V(H)$.
        Let $x$ be the neighbor of $z$ in $B$, and let $z' \in Z$ be the common neighbor of $z$ and $x$.            
        Then $x \in \phi(v)$.
    \end{claim}
    \begin{claimproof}
        Suppose not. Since $H$ is 2-connected, $v$ has at least two neighbors, say $u,w$. 
        Consequently, we must have $z' \in \phi(u)$ and $x \in \phi(w)$ (or vice versa, the other case is symmetric).
        As $H$ has at least four vertices, there is $q \in V(H) \setminus \{u,v,w\}$.
        However, removing $\phi(w)$ separates $\phi(v)$ from $\phi(q)$, so $w$ is a cutvertex in $H$, a contradiction.
    \end{claimproof}
    By the claim above we conclude that whenever some $z \in Z$ appears in a branch set of the induced minor model, 
    the same branch set contains the neighbor of $z$ in $B$.
    Consequently, we can modify the model into an equivalent one, by removing all vertices from $Z$ from all branch sets.
    Note that this preserves all (non-)adjacencies.
    Thus, $H$ is an induced minor of $K_{k,k} \odot K_2 - Z$, i.e., of $K_{k,k}$.
\end{proof}

\paragraph{Induced matching treewidth and forbidden substructures.}

Now let us analyze what substructures do \emph{not} appear in graphs listed in \cref{lem:constructions}.
We aim to detect maximal unbounded families, also keeping in mind the  comparison of containment (and exclusion) relations; recall \cref{fig:relations}.
Finally, whenever possible, we try to use constructions that are 2-connected, or at least of minimum degree 2.

\begin{lemma}\label[lemma]{lem:min-unbounded}
    The following classes of graphs have unbounded induced matching treewidth:
    \begin{enumerate}[(1.)]
        \item subcubic graphs of girth at least $t$, for any fixed $t$, \label{it:degree-girth}
        \item claw-free graphs, \label{it:clawfree}
        \item $H$-minor-free, where $H$ is not planar,\label{it:planar}        
        \item $P_4$-induced-minor free graphs, $(P_3+P_1)$-induced-minor-free graphs, $(K_3+P_1)$-induced-minor-free graphs,\label{it:join}
        \item bull-induced-subdivision-free graphs,\label{it:bull}
        \item house-induced-subdivision-free graphs,\label{it:house}
        \item $K_4$-induced-subdivision-free and $K'_4$-induced-minor-free graphs, where $K'_4$ is obtained from $K_4$ by subdividing one edge once,\label{it:k4}        
        \item $K_{2,2,2}$-induced-minor-free graphs, \label{it:K222}
        \item $H$-induced-minor-free graphs, where $H$ is any graph containing three vertices, each of degree at least 3, inducing $P_2+P_1$,  \label{it:P2P1-mindeg3}
        \item $H$-induced-minor-free graphs, where $H$ is any 2-connected graph containing $P_2 + P_1$,  \label{it:P2P1}
        \item 2-windmill-induced-minor-free graphs. \label{it:2windmill}
    \end{enumerate}
    Furthermore, items~\ref{it:degree-girth}--\ref{it:P2P1} hold even restricted to graphs of minimum at least degree 2,
    and items~\ref{it:degree-girth}--\ref{it:P2P1-mindeg3} hold even restricted to 2-connected graphs.
\end{lemma}
\begin{proof}
    We will use constructions provided in \cref{lem:constructions}.

    \paragraph{Walls and their line graphs.} Let us start with subdivided walls and their line graphs, as in \cref{lem:constructions}~\ref{it:wall} and~\ref{it:linewall}.

\smallskip
    \noindent \ref{it:degree-girth} Subdividing each edge of a large wall $t$ times we obtain a subcubic graph with girth at least $t$ and large induced matching treewidth.

\smallskip
     \noindent \ref{it:clawfree} Line graphs of subdivided walls are claw-free.

\smallskip
    \noindent \ref{it:planar} Note that every planar graph excludes any non-planar graph as a minor. Thus, (subdivided) walls are $H$-minor-free for any non-planar graph $H$.  

    \paragraph{Graphs $M_k$.} Fix $k \geq 2$ and consider the graph $M_k$ defined in \cref{lem:constructions}~\ref{it:Mk}. Let $X,Y$ be two sets that are complete to each other and induce $k$-edge matchings in $M_k$.

\smallskip
    \noindent \ref{it:join} It is straightforward to verify that for every $k\in\mathbb{N}$, the graph $M_k$ is $P_4$-free (i.e., a cograph) so, equivalently (see \cref{cor:linearforests}), $P_4$-induced-minor-free.
    Furthermore, as removing the neighborhood of any vertex leaves a matching, we observe that $M_k$ is $(P_3+P_1)$-induced-minor-free and $(K_3+P_1)$-induced-minor-free.

\smallskip
    \noindent \ref{it:bull} Suppose that for some $k$, the graph $M_k$ contains a copy of the bull as an induced subdivision.
    Denote the vertices of the bull by $a,b,c,u,v$, where $a,b,c$ form a triangle, $b$ is adjacent to $u$, and $c$ is adjacent to $v$.
    Since $u$ and $v$  are non-adjacent, we may assume by symmetry that their images $u',v'$ are both in $X$, and they belong to distinct matching edges in $M_k[X]$.
    Thus, the image $a'$ of $a$ must also be in $X$, belonging to yet another edge in $M_k[X]$.
    Since each of $u',v'$ is complete to $Y$, we observe that the images $b',c'$  of $b$ and $c$ are also in $X$.
    Among edges $ab$, $ac$ there is at least one, say $ab$, whose image is not contained in $X$ and thus, must be a three-vertex path $a'yb'$ for some $y \in Y$. However, $y$ is adjacent to $a'$ and $b'$, a contradiction.

    \noindent \ref{it:house} Suppose that for some $k$, the graph $M_k$ contains a copy of the house as an induced subdivision. Le $a,b,c,d$ be consecutive vertices of the 4-cycle in the bull, and let $z$ be the last vertex, adjacent to $a$ and $b$.
    As $M_k$ is $P_4$-free, it has no cycles with more than 4 vertices.
    Thus, the image of the cycle with vertices $a,b,c,d$ must be a 4-cycle in $M_k$.
    The only such cycles are obtained by taking two pairs of non-adjacent vertices, one from $X$ and the other from $Y$.
    However, there is no way to place the image of $z$ to make is non-adjacent both to the image of $c$ and $d$,
    a contradiction.

    \paragraph{Graphs $K_{k,k}^\circ$.} 
    Next, we will consider  the graph $K_{k,k}^\circ$ defined in \cref{lem:constructions}~\ref{it:Kkkcirc}. 
    Let $X$ and $Y$ denote, respectively, the sides of the copy of $K_{k,k}$ in $K_{k,k}^\circ$.

\smallskip
    \noindent  \ref{it:k4} 
    Note that for any $k\in\mathbb{N}$, the graph $K_{k,k}^\circ$ does not contain $K_4$. Thus, in order to show that it is $K_4$-inducedsubdivision-free, it is sufficient to show that it is $K'_4$-induced-minor-free.
    For contradiction, suppose that there is $k \geq 2$ such that $K_{k,k}^\circ$ contains $K'_4$ as an induced minor.
    Denote the degree-3 vertices of $K'_4$ by $a,b,c,d$, and let the degree-2 vertex $z$ be adjacent to $a$ and $b$.
    Consider an induced minor model $\phi : \{a,b,c,d,z\} \to 2^{V(K_{k,k}^\circ)}$.

    As $a$ has degree 3 in $K'_4$, we observe that $\phi(a)$ must intersect $X \cup Y$.
    By symmetry, suppose that $X \cap \phi(a) \neq \emptyset$.
    Similarly, we have that $\phi(b)$ must intersect $X \cup Y$.
    As $a$ and $b$ are non-adjacent, we observe that $\phi(a) \cap (X \cup Y) \subseteq X$ and $\phi(b) \cap (X \cup Y) \subseteq X$.
    Now, since $z$ is adjacent to both $a$ and $b$, we must have $\phi(z) \cap Y \neq \emptyset$.
    Consequently, $\phi(c)$ and $\phi(d)$ cannot intersect $X$.
    However, not there is no way to place $\phi(c)$ and $\phi(d)$ so that they are adjacent to each other, and also to $\phi(a)$ and $\phi(b)$, a contradiction.

\smallskip
   \noindent \ref{it:K222} For contradiction, suppose that there is $k\in\mathbb{N}$ such that $K_{k,k}^\circ$ contains $K_{2,2,2}$ as  an induced minor.
    Denote the vertices of $K_{2,2,2}$ by $\{a,a',b,b',c,c'\}$, where $aa',bb',cc'$ are non-edges.   
    Consider an induced minor model $\phi : \{a,a',b,b',c,c'\} \to 2^{V(K_{k,k}^\circ)}$.

    Let $u,u'$ be a pair of non-adjacent vertices of $K_{2,2,2}$.
    Note that it cannot hold that $\phi(u) \cap (X \cup Y) = \emptyset$, as $u$ is of degree 4 in $K_{2,2,2}$.
    By symmetry, suppose that $\phi(u) \cap X \neq \emptyset$.
    Note that since $u$ and $u'$ are non-adjacent in $K_{2,2,2}$, we must have $\phi(u') \cap (X \cup Y) \subseteq X$ and thus, symmetrically, $\phi(u) \cap (X \cup Y) \subseteq X$.

    Suppose by symmetry that $\phi(a) \cap (X \cup Y) \subseteq X$ and $\phi(a') \cap (X \cup Y) \subseteq X$.
    Note that since $b$ is adjacent to both $a$ and $a'$, we must have $\phi(b) \cap Y \neq \emptyset$,
    which means that $\phi(b) \cap (X \cup Y) \subseteq Y$ and $\phi(b') \cap (X \cup Y) \subseteq Y$.
    However, now there is no way to place $\phi(c)$ and $\phi(c')$, a contradiction.    

\smallskip
   \noindent \ref{it:P2P1-mindeg3} Now let $H$ be a graph as in \cref{it:P2P1-mindeg3},
    and denote the vertices of the induced $P_1 + P_2$ by $a,b,c$, where $bc$ is an edge. 
   For contradiction, suppose that there is $k\in\mathbb{N}$ such that $K_{k,k}^\circ$ contains $H$ as  an induced minor and consider an induced minor model $\phi: V(H) \to 2^{V(K_{k,k}^\circ)}$.

   As each of vertices $a,b,c$ is of degree at least 3, the image of each of them must intersect $X \cup Y$.
   By symmetry, assume that $\phi(a) \cap X \neq \emptyset$.
   As $a$ is not adjacent to both $b$ and $c$,
   it must hold that each of  $\phi(b) \cap (X \cup Y)$ and $\phi(c) \cap (X \cup Y)$ is a nonempty subset of $X$.   
   Since $b$ and $c$ are adjacent, there must be an edge between $\phi(a)$ and $\phi(a')$, 
   but this is impossible, a contradiction.

\paragraph{Graphs $K_{k,k} \odot K_p$ .} 
    In the next two items we will consider graphs $K_{k,k} \odot K_p$ for $p \in \{1,2\}$, defined in \cref{lem:constructions}~\ref{it:KkkodotK1} and~\ref{it:KkkodotK2}.  
    By $X,Y$ we denote the sides of the copy of $K_{k,k}$  in $K_{k,k} \odot K_p$, for $k \geq 2$ and $p \in \{1,2\}$.
   
\smallskip
\noindent \ref{it:P2P1} Let $H$ be a graph as in \cref{it:P2P1},
    and denote the vertices of the induced $P_1 + P_2$ by $a,b,c$, where $bc$ is an edge. 
    Consider $p=2$, i.e., the family of graphs $K_{k,k} \odot K_2$.
    Note that since $H$ is  2-connected and has at least 4 vertices,
    by \cref{lem:2connbiclique}, it is sufficient to show that $H$ is not an induced minor of $K_{k,k}$.
    So, for contradiction, suppose that there is an induced minor model $\phi : V(H) \to 2^{X \cup Y}$ of $H$ in $K_{k,k}$.

    By symmetry, suppose that there is $x \in X$ such that $x \in \phi(a)$.
    Since both $b$ and $c$ are non-adjacent to $a$, both sets $\phi(b)$ and $\phi(c)$ are contained in $X$.
    However, this means there are no edges between them, as $X$ is independent -- a contradiction with $bc$ being an edge.

\smallskip
    \noindent \ref{it:2windmill} We claim that $K_{k,k} \odot K_1$ is 2-windmill-induced-minor-free.
    Define $Z = V(K_{k,k} \odot K_1) \setminus (X \cup Y)$.
    Thus, $Z$ is a $2k$-element independent set and every vertex from $Z$ has a unique neighbor in $X \cup Y$.
    
    Consider a 2-windmill, and denote its central vertex by $c$, and the two edges obtained by removing $c$ by $aa'$ and $bb'$.
    For contradiction, suppose there is an induced minor model $\phi: \{c,a,a',b,b'\} \to 2^{X \cup Y \cup Z}$.

    Similarly as in \cref{clm:onlybiclique}, we observe that whenever some $z \in Z$ belongs to some $\phi(u)$ for $u \in \{c,a,a',b,b'\}$, then also $z'$, the unique neighbor of $z$, belongs to $\phi(u)$. Consequently, we can modify the minor model into an equivalent one by removing all vertices from $Z$ from all branch sets. Thus, the 2-windmill is an induced minor of $K_{k,k}$.
    In particular, $K_{k,k}$ contains $2K_2$ (the subgraph of the 2-windmill induced by $\{a,a',b,b'\}$) as an induced minor, so, by \cref{cor:linearforests}, also as an induced subgraph. This is a contradiction.        
\end{proof}

Finally, we state the following straighforward consequence of the transitivity of all considered containment relations. We will often use it implicitly.

\begin{observation}\label[observation]{obs:transfer-unboundedness}
    Let $\star$ be one of the considered relations, i.e., (induced) subgraph, (induced) topological minor, or (induced) minor.
    Let $F$ be contained in $H$ as $\star$, and let $\cG$ be a class of graphs.
    If induced matching treewidth is unbounded in the class of $F$-$\star$-free graphs from $\cG$,
    then it is unbounded in the class of $H$-$\star$-free graphs from $\cG$.
\end{observation}

\newpage
\section{Tipi-induced-minor-free graphs}\label{sec:tipi}
Let us start with the following structural observation concerning graphs that exclude 1-tipi as an induced minor.

\begin{lemma}\label[lemma]{lem:tipisubdiv}
    A graph $G$ is 1-tipi-induced-minor-free if and only if it is 1-tipi-induced-subdivision--free.
\end{lemma}
\begin{proof}    
    As every induced subdivision of any graph $H$ is in particular an induced minor model of $H$ (see \cref{fig:relations}), the forward implication is obvious.
    Thus suppose that $G$ has 1-tipi as an induced minor. We will exhibit an induced subdivision of 1-tipi in $G$.
    Denote the vertices of 1-tipi by $x,y,z$, and $p$, where $x,y,z$ form a triangle and $p$ is adjacent to $x$.
    Consider an induced minor model $\phi : \{x,y,z,p\} \to 2^{V(G)}$ of 1-tipi.
    By \cref{lem:minimal}, we can assume that $\phi(p)$ is a singleton, say, $\phi(p) = \{v\}$.
    
    Let $C$ be a shortest cycle that satisfies the following conditions:
    \begin{enumerate}
        \item $V(C)$ consists of three consecutive non-empty subpaths, each contained in, respectively, $\phi(x)$,$\phi(y)$, and $\phi(z)$,        
        \item $V(C)$ contains a neighbor of $v$.
    \end{enumerate}
    Such a cycle exists since $\phi$ is an induced minor model of 1-tipi.
    For convenience, fix the ordering of the cycle from $\phi(y)$ to $\phi(x)$, and then to $\phi(z)$.
    
    First, suppose that $v$ is adjacent to at least two vertices of $C$.
    Let $a,b$ be two consecutive vertices of $C$ such that $v$ is non-adjacent to $a$ and adjacent to $b$,
    and $b$ is later than $a$ along the cycle.
    Such vertices exist as $v$ has a neighbor in $\phi(x)$ and has no neighbor in $\phi(y)$.
    Let $c$ be the first neighbor of $v$ that appears on $C$ after $b$.
    Note that $ac \notin E(G)$, as the neighborhood of $v$ in $C$ is contained in $\phi(x)$ and $C$ intersects both $\phi(y)$ and $\phi(z)$.
    Now, vertices $a,b,c,v$ and the subpath of $C$ joining $b$ and $c$ (and avoiding $a$) form an induced subdivision of 1-tipi.

    So suppose that $v$ has exactly one neighbor on $C$, say $u$.
    Since $C$ is shortest possible, we observe that either $C$ is an induced cycle or there is an edge between the neighbors of $u$ on $C$.
    In both outcomes we obtain an induced subdivision of 1-tipi.    
\end{proof}

We also give the following easy observation that will be used in this section but also in the proof of~\cref{thm:windmill}.

\begin{observation}\label[observation]{obs:easywindmill}
    Let $G$ be a graph that contains a matching $M$ of size $k\ge 1$ as an induced subgraph. Also, let $Z$ be a connected subset of $V(G)\setminus V(M)$. If every vertex of $V(M)$ is adjacent to some vertex of $Z$, then $G$ contains the $k$-windmill as an induced minor. In particular, $G$ contains the $(k-1)$-tipi as an induced minor.
\end{observation}

\begin{lemma}\label[lemma]{lemma:blackbox-tipi}
    Let $k\ge 1$, let $G$ be a graph, and let $Z\subseteq V(G)$ such that $G[Z]$ is connected and contains a cycle $C$ as a subgraph. Also, suppose that $G$ contains a matching $M$ of at least $3k$ edges as an induced subgraph so that for every edge $e$ of $M$ one of the following holds:
    \begin{itemize}
        \item both endpoints of $e$ are disjoint from $Z$ and at least one of them is adjacent to a vertex in $Z$, or
        \item $e$ has exactly one endpoint in  $Z$, which is also in $V(C)$.
    \end{itemize}
    Then, $G$ contains the $k$-tipi as an induced minor.
\end{lemma}
\begin{proof}
    An edge of $M$ can satisfy exactly one of the conditions in the statement. Let $E_1$ be the edges in $M$ that satisfy the first condition and let $E_2$ be the edges that satisfy the second one.

    First, suppose there is an edge $e \in E_1$ whose both endpoints are adjacent to vertices in $Z$.
    In that case,  from each edge in $E_1 \setminus \{e\}$ we delete one endpoint, leaving one that has a neighbor in $Z$, and then we contract $Z$ to a single vertex.
    This way we obtain a $(|M|-1)$-tipi (and thus, a $k$-tipi) as an induced minor of $G$.

    So from now on, assume that every edge $e\in E_1$ has exactly one endpoint that is adjacent to a vertex in $Z$. We denote such endpoint by $u_e$.
    
    For each edge $e\in E_1$, we fix a shortest path $P_e$ in $G$ joining $u_e$ to some vertex in $V(C)$.
    Note that $V(P_e)\setminus \{u_e\}\subseteq Z$.
    By $w_e$ we denote the endpoint of $P_e$ in $V(C)$.
    For edges $e\in E_2$, we set $w_e$ be their endpoint in $C$.
    Note that the paths $P_e$ do not contain any edge from $E_2$. 

    Next, we fix two consecutive vertices $z_1$ and $z_2$ in $C$ and we partition $V(C)$ to the sets $B_1 = \{z_1\}$, $B_2=\{z_2\}$, and $B_3=V(C)\setminus \{z_1,z_2\}$.
    Fix $j\in\{1,2,3\}$ so that the number of edges $e\in M$ with $w_e\in B_j$ is maximized.
    Since $|M|\ge 3k$, $B_j$ contains the vertices $w_e$ of at least $k$ edges in $M$.
    Thus, by contracting the paths $P_e$ for each $e\in E_1$ with $w_e\in B_j$ to single edges, contracting $V(C) \setminus \{z_1,z_2\}$ to a single vertex, and removing all remaining vertices, we obtain the $k'$-tipi for $k' \geq k$ (and thus, also the $k$-tipi) as an induced minor.
\end{proof}

As the main result of this section, we show excluding a tipi as an induced minor yields a class of bounded induced matching treewidth.

\begin{theorem}\label{thm:tipi}
    For every integer $k\ge 1$, every graph $G$ that does not contain the $k$-tipi as an induced minor satisfies $\yolov(G) \leq 11k+5$.
\end{theorem}
\begin{proof}
    \renewcommand{\claimprefix}{\ref{thm:tipi}}
    \setcounter{claim}{0}

    Without loss of generality, assume that $G$ is connected, as otherwise we can process each component separately and otain a tree decomposition of $G$ by combining the tree decompositions of components.
    If $G$ is a tree, then $\tw(G) \leq 1$ and thus, due \cref{prop:properties}~\ref{it:yolov-tin-tw}, we have $\yolov(G)\le 2$.
    So, assume hat $G$ has a cycle.
    
    Pick any vertex $r$ of $G$ and consider the BFS-layering  $(L_0,\ldots,L_q)$ of $G$ rooted at this vertex.
    
    \begin{claim}
        For every $i \ge 0$, every induced matching in $G[L_i]$ has at most $k$ edges.
    \end{claim}
    \begin{claimproof}
        The statement is trivial for $i=0$, so assume that $i \geq 1$.
        For contradiction, suppose there is an induced matching $M$ on $k+1$ edges contained in $L_i$.
        Note that the set $Z=\bigcup_{j<i}L_j$ is connected and disjoint from $L_i$.
        Furthermore, every vertex in $L_i$ has a neighbor in $Z=\bigcup_{j<i}L_j$.        
        Thus, by~\cref{obs:easywindmill}, $G$ contains the $k$-tipi as an induced minor, a contradiction.
    \end{claimproof}

    We next bound the size of an induced matching between layers.
    
    \begin{claim}
        For every $i\in[q]$, every induced $(L_{i-1},L_{i})$-matching of $G$ has at most $3k+1$ edges.
    \end{claim}
    \begin{claimproof}
    We assume that $i\ge 2$ since in the case $i=1$ we have $|L_0|=1$.
    Suppose that $G$ contains an induced $(L_{i-1},L_{i})$-matching $M$ of at least $3k+2$ edges.
    Recall that we assumed that $G$ contains a cycle as a subgraph. Let $C$ be a shortest such cycle; in particular, it is induced.
    We distinguish two cases depending whether $C$ uses some edge from $M$.
    In both cases, we use~\cref{lemma:blackbox-tipi} to arrive at a contradiction.
    
    \paragraph{Case A: $E(C)\cap M=\emptyset$.}  
    Let $j_0=\min\{j\in[q] ~|~ V(C)\cap L_j\neq\emptyset\}$ and let $z\in V(C)\cap L_{j_0}$. Let $P_z$ be a root path for $z$.
    Note that there is at most one edge $e'\in M$ that has an endpoint in $V(P_z)$, as $P_z$ contains at most one vertex from $L_{i-1}$ and at most one vertex from $L_i$, and these two vertices are adjacent.
    Let $M' = M \setminus \{e'\}$ (or $M' = M$ is $e'$ does not exist), note that $|M'| \geq 3k$.
    For every edge in $M'$ we have that both endpoints are disjoint from $V(P_z)$.

    Let $Z=\bigcup_{j<i-1}L_{j}\cup V(C)\cup V(P_z)$ and note that $G[Z]$ is connected and contains the cycle $C$. Our goal is to apply~\cref{lemma:blackbox-tipi} for $Z$ and $M'$.
    Since $E(C)\cap M=\emptyset$ and $C$ is induced, every edge in $M'$ has at most one endpoint in $V(C)$. Let $E_1=\{uv\in M' ~|~ \{u,v\}\cap V(C)=\emptyset\}$ and $E_2=M'\setminus E_1$.
    Note that for every edge $e\in E_1$, both endpoints are disjoint from $Z$ and the endpoint in $L_{i-1}$ is adjacent to $Z$.
    For every edge in $E_2$, only one of its endpoints is in $V(C)$ and none is in $Z\setminus V(C)$. Therefore, by applying~\cref{lemma:blackbox-tipi} for $Z$ and $M'$, we get that $G$ contains the $k$-tipi as an induced minor, a contradicion.

    \paragraph{Case B: $E(C)\cap M\neq\emptyset$.} Note that in this case there are at least two edges between vertices in $L_{i-1}$ and $L_{i}$ that belong to $E(C)$.

     Consider an inclusion-wise minimal subpath $Q$ of $C\cap \bigcup_{j>i-1}L_j$ with both endpoints in $L_{i}$ and all internal vertices in $\bigcup_{j>i}L_j$.
     Let $v_1,v_2$ be the endpoints of $Q$ and let $u_1$ and $u_2$ be the neighbors of $v_1$ and $v_2$ respectively in $C$ which belong to $L_{i-1}$; if $v_1=v_2$, pick $u_1$ and $u_2$ be distinct vertices.
     Then, let $P_1$ and $P_2$ be, respectively, root paths of $u_1$ and $u_2$. Now set $Z=\bigcup_{j < i-1} L_j \cup V(Q)\cup V(P_1)\cup V(P_2)$ and note that $G[Z]$ is connected and contains a cycle.
     Also note that by construction there are at most two edges in $M$ that have some endpoint in $Z$. 
     Let $M'$ be obtained from $M$ by removing these two edges (if they exist); note that $|M'| \geq 3k$.
     Furthermore, for every edge from $M'$, its endpoint in $L_{i-1}$ has a neighbor in $L_{i-2} \subseteq Z$.
     Therefore, by applying~\cref{lemma:blackbox-tipi} to $Z$ and $M'$ we get that $G$ contains a $k$-tipi as an induced minor, a contradiction.
    \end{claimproof}

    Thus, by \cref{lem:from-bfs-to-tree} for $a=k$ and $b = 3k+1$, we obtain that $\yolov(G) \leq 2k + 3(3k+1) = 11k+3$. This completes the proof.
\end{proof}

\newpage
\section{Dichotomies for one excluded graph}\label{sec:dichos}
In this section we prove dichotomies for boundedness of induced matching treewidth in classes that exclude a fixed substructure.
For exclusion as a subgraph, as an induced subgraph, as a minor, or as a topological minor, the dichotomies almost coincide with their counterparts for tree-independence number.

\subsection{... as a subgraph or an induced subgraph}

    \begin{theorem}[\cref{thm:dichoeasy}~\ref{it:easy-subgraph} restated]
        The class of $H$-subgraph-free graphs has bounded induced matching treewidth if and only if every component of $H$ is a subcubic tree with at most one vertex of degree 3.        
        All unbounded cases remain so, even when restricted to 2-connected graphs.
\end{theorem}
    \begin{proof}   
        If every component of $H$ is a subcubic tree with at most one vertex of degree 3, then $H$-subgraph-free graphs have bounded tree-independence number (see \cref{tab:summary}), and thus also induced matching treewidth.

        Otherwise $H$ contains a cycle or a vertex of degree at least 4.
        In this case, subcubic graphs of girth at least $|V(H)|+1$ are $H$-subgraph-free and have unbounded induced matching treewidth by \cref{lem:min-unbounded}~\ref{it:degree-girth}.        
\end{proof}

    \begin{theorem}[\cref{thm:dichoeasy}~\ref{it:easy-inducedsubgraph} restated]
        The class of $H$-free graphs has bounded induced matching treewidth if and only if $H = P_3$ or $H$ has maximum degree at most 1 (i.e., is a subgraph of a matching). 
        All unbounded cases remain so, even when restricted to 2-connected graphs.
\end{theorem}
    \begin{proof}   
        Recall that $P_3$-free graphs are precisely graphs whose every component is a clique; clearly they have bounded induced matching treewidth.
        If $H$ is a subgraph of a matching, then $H$-free graphs do not contain large induced matchings at all, so in particular induced matching treewidth is bounded.
        
        Now, let us assume that $H$ is neither $P_3$ nor a graph with maximum degree at most 1.
        Each such $H$ contains one of the following graphs as an induced subgraph: $K_3$, $P_3+P_1$, $P_4$, $K_{1,3}$, $C_4$.
        Each of them can be excluded by \cref{lem:min-unbounded}~\ref{it:degree-girth},~\ref{it:clawfree}, and~\ref{it:join}.
\end{proof}

\subsection{... as a minor or a topological minor}
\begin{theorem}[\cref{thm:dichoeasy}~\ref{it:easy-minor} restated]
        The class of $H$-minor-free graphs has bounded induced matching treewidth if and only if $H$ is planar.
        All unbounded cases remain so, even when restricted to 2-connected graphs.
\end{theorem}
\begin{proof}
    If $H$ is planar, then $H$-minor-free graphs have bounded treewidth~\cite{RobertsonS86} and thus also induced matching treewidth.
    Otherwise, induced matching treewidth is unbounded by~\cref{lem:min-unbounded}~\ref{it:planar}.
\end{proof}

\begin{theorem}[\cref{thm:dichoeasy}~\ref{it:easy-topminor} restated]
    The class of $H$-topological-minor-free graphs has bounded induced matching treewidth if and only if $H$ is subcubic and planar.    
     All unbounded cases remain so, even when restricted to 2-connected graphs.
\end{theorem}
\begin{proof}
    It is well-known that if $H$ is subcubic, then classes of $H$-minor-free graphs and $H$-topological-minor-free graphs coincide (see, e.g., Diestel~\cite{DBLP:books/daglib/0030488-6}). 
    Thus, if $H$ is subcubic and planar, then $H$-topological-minor-free (equivalently, $H$-minor-free graphs) graphs have bounded treewidth~\cite{RobertsonS86} and thus also induced matching treewidth.
    Otherwise, induced matching treewidth is unbounded by~\cref{lem:min-unbounded}~\ref{it:degree-girth} and~\ref{it:planar}.
   
\end{proof}

\subsection{... as an induced subdivision}
Now let us consider $H$-induced-subdivision-free graphs, for any fixed $H$.

\thmdichosubdivision*
\begin{proof}
    \renewcommand{\claimprefix}{\ref{thm:dichosubdivision}}
    \setcounter{claim}{0}

    If $H$ is an induced subdivision of $C_4$ or $K_4^-$, then $H$-induced-subdivision-free graphs have bounded tree-independence number (see~\cref{tab:summary}), and thus also induced matching treewidth.
    If $H$ is 1-tipi, then the claim follows from combining \cref{lem:tipisubdiv} and \cref{thm:tipi}.
    If $H$ is a graph of maximum degree 1, then, by \cref{cor:linearforests},
    $H$-induced-subdivision-free graphs coincide with $H$-free graphs.
    Such graphs do not have a large induced matching at all, so in particular they have bounded induced matching treewidth.

    \bigskip
    Thus, let us argue that all remaining cases define a class with unbounded induced matching treewidth.
    We employ a small case analysis. We remark that in all subsequent claims we assume that $H$ satisfies all properties established so far.
    Let us start with a few simple observations.

    \begin{claim}\label{inducedsubdiv-easy}
        Induced matching treewidth of $H$-induced-subdivision-free graphs is unbounded if:
        \begin{itemize}
            \item $H$ contains $K_{1,3}$ as an induced subgraph,
            \item $H$ contains $P_3 + P_1$ or $P_4$ or $K_3 + P_1$ as an induced subdivision,
            \item $H$ is non-planar,
            \item $H$ has a vertex of degree at least 4.
        \end{itemize}
    \end{claim}
    \begin{claimproof}    
    First, if $H$ contains $K_{1,3}$ as an induced subgraph, then claw-free graphs are $H$-induced-subdivision-free and thus induced matching treewidth is unbounded by \cref{lem:min-unbounded}~\ref{it:clawfree}.
    Second, if $H$ contains $(P_3 + P_1)$ or $P_4$ or $(K_3+P_1)$ as an induced subdivision, then it also contains one of these graphs as an induced minor, so \cref{lem:min-unbounded}~\ref{it:join} yields that induced matching treewidth is unbounded in $H$-induced-subdivision-free graphs.
    Third, every $H$-minor-free graph is $H$-topological-minor-free, and every $H$-topological-minor-free graph is $H$-induced-subdivision-free.
    If $H$ is non-planar, then the former class has unbounded induced matching treewidth by \cref{lem:min-unbounded}~\ref{it:planar}.
    If $H$ has a vertex of degree at least 4, then every subcubic graph is $H$-topological-minor-free, so the latter class has unbounded induced matching treewidth by \cref{lem:min-unbounded}~\ref{it:degree-girth}.
    \end{claimproof}
   
    So suppose that $H$ is a subcubic planar graph that does not contain $K_{1,3}$ as an induced subgraph nor $P_4$, $P_3 + P_1$, $K_3+P_1$ as an induced subdivision.
    If $H$ is a forest, then these conditions imply that $H$ is either $P_3$ (and thus, contained in $C_4$ as an induced subdivision) or has maximum degree 1 -- in both cases we have one of bounded cases considered above. 
    Thus, $H$ is not a forest. Furthermore, since $H$ excludes $K_3 + P_1$ as an induced subdivision, every cycle in $H$ is dominating.

    \begin{claim}
        If $H$ has an induced cycle on at least four vertices, then induced matching treewidth of $H$-induced-subdivision-free graphs is unbounded.
    \end{claim}
    \begin{claimproof}
    If $H$ has an induced cycle on at least 5 vertices, then it also contains $P_3+P_1$, so again we have unbounded case.
    Suppose that $H$ has an induced cycle on four vertices $\{a,b,c,d\}$.
    Since $H$ is not isomorphic to $C_4$, there must be a vertex $x \in V(H) \setminus \{a,b,c,d\}$.
    Since every cycle is dominating, $x$ must have neighbors in $\{a,b,c,d\}$.
    If $x$ has one neighbor in $\{a,b,c,d\}$, then $H$ contains $P_3+P_1$, a contradiction.
    Similarly, if $x$ has two neighbors on $\{a,b,c,d\}$ and they are non-adjacent, then $H$ contains an induced $K_{1,3}$.

    If the neighborhood of $x$ in $\{a,b,c,d\}$ consists of two adjacent vertices, then $H$ contains a house as an induced subdivision, and thus induced matching treewidth is unbounded by \cref{lem:min-unbounded}~\ref{it:house}.
    If the neighborhood of $x$ in $\{a,b,c,d\}$ consists of three vertices, then $H$ contains $K_4$ as an induced subdivision and induced matching treewidth is unbounded by \cref{lem:min-unbounded}~\ref{it:k4}.

    As $H$ is subcubic, this covers all cases if $H$ contains an induced $C_4$.
    \end{claimproof}

    Thus, assume that $H$ is not a forest and all induced cycles in $H$ are triangles, i.e., $H$ is chordal.
    If $H$ contains $K_4$ as an induced subdivision, then $H$-induced-subdivision-free graphs have unbounded induced matching treewidth by \cref{lem:min-unbounded}~\ref{it:k4}.
    We observe that every such graph that does not contain previously excluded substructures contains either a $K_4^-$ or a 1-tipi (as an induced subgraph).
    We know that $H$ must contain at least one more vertex $v$, as both these structures were considered in the first part of the proof and now we assume that $H$ is not one of such graphs.

    \begin{claim}
        If $H$ contains  a $K_4^-$, then induced matching treewidth of $H$-induced-subdivision-free graphs is unbounded. 
    \end{claim}
    \begin{claimproof}
    Fix one copy of $K_4^-$ in $H$.
    Since $H$ is $(K_3+P_1)$-free, $v$ is adjacent to both triangles in $K_4^-$.
    Since $H$ is subcubic, the neighborhood of $v$ in the considered copy of $K_4^-$ consists of both degree-2 vertices.
    However, this means that $H$ contains an induced $C_4$, a contradiction.
    \end{claimproof}

    \begin{claim}
        If $H$ contains  a 1-tipi, then induced matching treewidth of $H$-induced-subdivision-free graphs is unbounded. 
    \end{claim} 
    \begin{claimproof}
    Fix one copy of 1-tipi in $H$.
    Since $H$ is $(K_3+P_1)$-free and subcubic, $v$ is adjacent to at least one of degree-2 vertices of the 1-tipi.
    If $v$ is also adjacent to the degree-1 vertex of the 1-tipi, $H$ has an induced $C_4$, a contradiction.
    If $v$ is adjacent to both degree-2 vertices of the 1-tipi, then $H$ is not $(K_3+P_1)$-free, a contradiction.
    Thus, the neighborhood of $v$ in the considered copy of 1-tipi consists of exactly one of degree-2 vertices. Then the vertices of the copy of 1-tipi, together with $v$, induce in $H$ a copy of the bull. and we get an unbounded case by \cref{lem:min-unbounded}~\ref{it:bull}.
    \end{claimproof}
    
    Summing up, for every graph $H$ not listed in the first part of the statement, induced matching treewidth of $H$-induced-subdivision-free graphs is unbounded.
    As we only used constructions from \cref{lem:min-unbounded}~\ref{it:degree-girth}--\ref{it:K222}, the statement holds even when we restrict attention to 2-connected graphs.
    This completes the proof.
\end{proof}

\subsection{... as an induced minor}
Finally, we consider the last contaiment relation, i.e., induced minors.

\thmdichoindminor*
\begin{proof}
    \renewcommand{\claimprefix}{\ref{thm:dichoindminor}}
    \setcounter{claim}{0}

    If $H$ is an induced minor of $K_{2,q}$ for some fixed $q$ or $W_4$, or $K_5^-$, then $H$-induced-minor-free graphs have bounded tree-independence number (see~\cref{tab:summary}), and thus also induced matching treewidth.
    If $H$ is a tipi, then the claim follows from \cref{thm:tipi}.
    If $H$ is a graph of maximum degree 1, then, by \cref{cor:linearforests}, classes of $H$-free graphs and $H$-induced-minor-free graphs coincide and such graphs do not have large induced matchings at all.

    \bigskip
    Now we need to show that all remaining graphs $H$ define classes of unbounded induced matching treewidth.
    Thus, assume that $H$ is not one of the graphs listed above.

    \begin{claim}
    If $H$ is disconnected, then induced matching treewidth of  $H$-induced-minor-free graphs is unbounded.
    \end{claim}
    \begin{claimproof}
    Suppose $H$ is disconnected.
    Since $H$ is not a subgraph of a matching, at least one of components of $H$, say $C$, must have at least three vertices.
    Thus, $C$ contains either $P_3$ or $K_3$ as an induced subgraph.
    Consequently, $H$ contains either $P_3 + P_1$ or $K_3 + P_1$ as an induced minor, so the conclusion follows from \cref{lem:min-unbounded}~\ref{it:join}.
    \end{claimproof}

    Next, we consider two cases, depending whether $H$ is 2-connected or not.

\paragraph{Case A. $H$ is not 2-connected.}
    Let $\cB$ and $\cC$ denote, respectively, the set of blocks and the set of cutvertices of $G$.
    Let $T$ be the block-cut tree of $G$ .   

    \begin{claim}
    If $H$ has at least two cutvertices, then induced matching treewidth of  $H$-induced-minor-free graphs is unbounded.
    \end{claim}
    \begin{claimproof}        
        If $H$ has at least two cutvertices, $T$ contains a 5-vertex path with three vertices from $\cB$ (and two from $\cC$) -- recall that every leaf of $T$ is in $B$.
        Denote the consecutive vertices of this path by $x,a,y,b,z$.
        Let $v_x$ (resp., $v_z$) be a neighbor of $a$ (resp., $b$)  in the block $x$ (resp., $z$).
        Let $P$ be a shortest $a$-$b$-path contained in $y$.
        Note that $v_x,a,P,b,v_z$ is an induced path with at least four vertices.
        Thus, $P_4$ is an induced minor of $H$, and the conclusion follows from \cref{lem:min-unbounded}~\ref{it:join}.
    \end{claimproof}

    So, assume that $H$ has exactly one cutvertex, denote it by $c$.
    Consequently, $T$ is a star with $c$ as the central vertex.    
    
    \begin{claim}
        If $H$ has a block with at least four vertices, then induced matching treewidth of $H$-induced-minor-free graphs is unbounded.
    \end{claim}
    \begin{claimproof}
        Let $b$ be a  block of $G$ with at least four vertices and let $b'$ be some block other than $b$;
        it exists since $G$ is not 2-connected.
        As $b - c$ is a connected graph on at least three vertices, it contains an induced $P_3$ or $K_3$.
        This subgraph, together with any vertex from $b' - c$, forms an induced $P_3 + P_1$ or induced $K_3 + P_1$ in $H$.
        In both cases, the conclusion follows from \cref{lem:min-unbounded}~\ref{it:join}.
    \end{claimproof}

    Thus, every block of $G$ has either two or three vertices, i.e., is either a $K_2$ or $K_3$.
    
\begin{claim}\label{clm:2windmill}
        If $H$ has two blocks with three vertices, then induced matching treewidth of $H$-induced-minor-free graphs is unbounded.
    \end{claim}
    \begin{claimproof}
        If $H$ has two blocks with three vertices, then it contains a 2-windmill as an induced subgraph.
        The conclusion follows from \cref{lem:min-unbounded}~\ref{it:2windmill}.
    \end{claimproof}

    Consequently, $H$ has at most one block with three vertices.
    This means that $H$ is an induced subgraph of a tipi -- a graph whose exclusion defines a class of bounded induced matching treewidth.
    This completes the proof in Case A.

\paragraph{Case B. $H$ is 2-connected.}
    First, we observe the following.

    \begin{claim}\label{clm:p2p1}
        If $H$ contains induced $P_2 + P_1$, then induced matching treewidth of $H$-induced-minor-free graphs is unbounded.
    \end{claim}
    \begin{claimproof}
        The claim follows from \cref{lem:min-unbounded}~\ref{it:P2P1}.
    \end{claimproof}

    So from now on let us assume that $H$ is $(P_2+P_1)$-free, i.e., is a complete multipartite graph.
    Let $k$ be the number of parts (i.e., inclusion-wise maximal independent sets) in $H$.
    A part is \emph{large} if it has at least two vertices.

    \begin{claim}\label{clm:threelarge}
        If $H$ has at least 3 large parts, then induced matching treewidth of $H$-induced-minor-free graphs is unbounded.
    \end{claim}
    \begin{claimproof}
        If $H$ has at least three large parts, then it contains $K_{2,2,2}$ as an induced subgraph.
        The claim follows from \cref{lem:min-unbounded}~\ref{it:K222}.
    \end{claimproof}
     
    \begin{claim}
        If $k \geq 5$, then induced matching treewidth of $H$-induced-minor-free graphs is unbounded.
    \end{claim}
    \begin{claimproof}
        If $k \geq 5$, then $H$ contains $K_5$ and thus, is non-planar.
        The claim  follows from \cref{lem:min-unbounded}~\ref{it:planar}.
    \end{claimproof}

    So, $H$ has two, three, or four parts. Let us consider these cases separately.

    \begin{claim}
     If $k=2$, then either $H$ is an induced minor of $K_{2,q}$ for some $q$, or induced matching treewidth of $H$-induced-minor-free graphs is unbounded.
    \end{claim}
    \begin{claimproof}
        If each part has at least three vertices, then $H$ contains $K_{3,3}$ and thus, is non-planar, so we are done by \cref{lem:min-unbounded}~\ref{it:planar}.
        Otherwise $H$ is an induced minor of $K_{2,q}$ for some $q \geq 2$, which is one of the cases that imply bounded induced matching treewidth.
    \end{claimproof}

    \begin{claim}
        If $k=3$, then either $H$ is an induced minor of $K_{2,q}$ for some $q$ or isomorphic to $W_4$, or induced matching treewidth of $H$-induced-minor-free graphs is unbounded.
    \end{claim}
    \begin{claimproof}
        By \cref{clm:threelarge} we know that at least one part is just a singleton.        
        Let $H = K_{a,b,1}$, where $a \geq b$.
        
        If $b = 1$, then $H$ is an induced minor of $K_{a+1,2}$,
        i.e., is one of the graphs that defines a class of bounded induced matching treewidth.
        If $a \geq 3$, then $H$ contains $K_{3,3}$ as a subgraph and thus is non-planar, so we are done by \cref{lem:min-unbounded}~\ref{it:planar}.
        So, we are left with the case that $H = K_{2,2,1} \simeq W_4$.
        In this case, we obtain a class of bounded induced matching treewidth.
    \end{claimproof}

    \begin{claim}
        If $k=4$, then either $H$ is isomorphic to $K_5^-$, or induced matching treewidth of $H$-induced-minor-free graphs is unbounded.
    \end{claim}
    \begin{claimproof}
        By \cref{clm:threelarge} we know that two parts are singletons.        
        Let $H = K_{a,b,1,1}$, where $a \geq b$.
        If $a \geq 3$ or $a=b=2$, then $H$ contains $K_{3,3,}$ as a subgraph and again we are done since $H$ is non-planar.
        So, we are left with the case that $H = K_{2,1,1,1} \simeq K_5^-$.
        In this case, we obtain a class of bounded induced matching treewidth.
    \end{claimproof}
    This completes the proof.    
\end{proof}

Let us remark that the case that  \cref{clm:2windmill} and \cref{clm:p2p1} are the only cases in the proof of \cref{thm:dichoindminor} where we use constructions that are not 2-connected, i.e., \cref{lem:min-unbounded}~\ref{it:P2P1} and~\ref{it:2windmill}.
Furthermore, \cref{clm:2windmill} is the only place  where we use \cref{lem:min-unbounded}~\ref{it:2windmill}, i.e., graphs that contain vertices of degree 1.
This gives the following observation.

\begin{observation}\label[observation]{obs:afterindminor}
Let $H$ be a graph for which $H$-induced-minor-free graphs have unbounded induced matching treewidth.
\begin{enumerate}
    \item Induced matching treewidth in $H$-induced-minor-free graphs remains unbounded, even restricted to graphs of minimum degree at least 2, unless $H$ contains a 2-windmill as an induced subgraph.
    \item Induced matching treewidth in $H$-induced-minor-free graphs remains unbounded, even restricted to 2-connected graphs, unless $H$ contains a 2-windmill as an induced subgraph or $H$ is 2-connected and contains $P_2 + P_1$ as an induced subgraph.
\end{enumerate}
\end{observation}

\cref{obs:afterindminor} will be a starting point to further discussion in \cref{sec:noleaves} and \cref{sec:2connected}.

\newpage
\section{Restricting to graphs of minimum degree at least 2}\label{sec:noleaves}
Now, let us discuss the (un)boundedness of the considered parameters when we additionally assume that graphs have no vertices od degree 1.
By \cref{lem:block-reduction}, this restriction does not make any difference for tree-independence number.
Furhermore, recall that this is also the case for induced matching treewidth and graphs that excluding a fixed graph $H$ as an (induced) subgraph, (induced) topological minor, or a minors.
However, it turns our that for $H$-induced-minor-free graphs, the situation changes.

\subsection{Windmill-induced-minor-free graphs of minimum degree 2}

We show that if we assume that no vertices have degree at most 1, then we get one more family of graphs that define classes of bounded induced matching treewidth.

\begin{theorem}\label{thm:windmill}
    There is a function $f:\mathbb{N}\to\mathbb{N}$ such that for every integer $k\ge 0$, every graph $G$ of minimum degree at least 2 that does not contain the $k$-windmill as an induced minor  satisfies $\yolov(G) \leq 36k^3+2k$.
\end{theorem}
\begin{proof}
    \renewcommand{\claimprefix}{\ref{thm:windmill}}
    \setcounter{claim}{0}

    Without loss of generality we may assume that $G$ is connected.    
    As in the proof of~\cref{thm:tipi}, consider the BFS-layering $(L_0,\ldots,L_q)$ of $G$ rooted at an arbitrary vertex $r$. We aim to bound the size of a largest induced matching within each layer and between any two consecutive layers, and then apply \cref{lem:from-bfs-to-tree}.
    
    \begin{claim}
       For every $i \geq 0$, every induced matching of $G$ between vertices of $L_i$ has at most $k-1$ edges.
    \end{claim}
    \begin{claimproof}
        For $i=0$ the statement is trivial, and for $i \geq 1$, as in the proof of~\cref{thm:tipi},
        the claim follows by applying~\cref{obs:easywindmill} for $Z=\bigcup_{j<i}L_j$.
    \end{claimproof}

    We next bound the size of an induced matching between consecutive layers.
    
    \begin{claim}
        For every $i \in [q]$, every induced $(L_{i-1},L_i)$-matching of $G$ has at most $12k^3-1$ edges.
    \end{claim}
    \begin{claimproof}
    We assume that $i\ge 2$ since in the case $i=1$ we have $|L_0|=1$.
    We fix $i\in[2,q]$ and suppose that $G$ contains an induced $(L_{i-1},L_i)$-matching $M$ with at least $12k^3$ edges.
    We denote by $M^{\textrm{up}}$ (resp. $M^{\textrm{down}}$) the set of endpoints of $M$ in $L_{i-1}$ (resp. $L_i$) and by $G'$ the graph $G[\bigcup_{j \geq i} L_{j}]$ -- see~\cref{fig:MupMdown} for an illustration.
    \begin{figure}[ht]
         \centering
         \begin{tikzpicture}[every node/.style={minimum size = 4pt, inner sep = 0pt}]
            \node[circle,draw=black,fill=white] (left1) at (0,0) {};
            \node[circle,draw=black,fill=white] (left2) at (0,-1) {};
            \draw[blue,line width=1pt] (left1) to (left2);
            \node[circle,draw=black,fill=white] (c1) at (1,0) {};
            \node[circle,draw=black,fill=white] (c2) at (1,-1) {};
            \draw[blue,line width=1pt] (c1) to (c2);
            \node[circle,draw=black,fill=white] (right1) at (2,0) {};
            \node[circle,draw=black,fill=white] (right2) at (2,-1) {};
            \draw[blue,line width=1pt] (right1) to (right2);
            \node[label={right:$M^{\textrm{up}}$}] at (2.5,0) {};
            \node[label={right:$M^{\textrm{down}}$}] at (2.5,-1) {};

            \draw[] (-1,-2) .. controls (-0.5,-0.2) and (2.5,-0.2) .. (3,-2);
            \node[label={center:$G'$}] at (0,-1.5) {};
    \end{tikzpicture}
    \caption{An illustration of the sets $M^{\textrm{up}}$ and $M^{\textrm{down}}$ and the graph $G'$.}
    \label{fig:MupMdown}
    \end{figure}
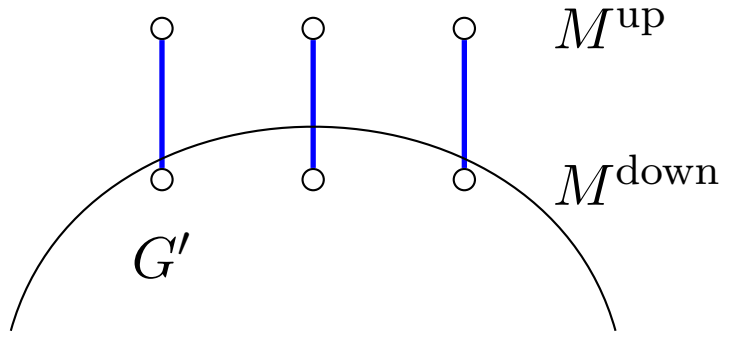
    
    We distinguish cases:

    \paragraph{Case A: There is a component $C$ of $G'$ containing at least $4k^2$ vertices of $M^{\textrm{down}}$.}

    In this case, consider an inclusion-wise minimal (with respect to the vertex set) tree $T$ in $C$ which satisfies the following properties:
    \begin{itemize}
        \item $T$ spans $M^{\textrm{down}}\cap V(C)$ and
        \item among all inclusion-wise minimal trees that satisfy the property above, $T$ maximizes the number of leaves that are in $M^{\textrm{down}}\cap V(C)$.
    \end{itemize}
    Note that every leaf of $T$ is in $M^{\textrm{down}}$ and therefore the leaves of $T$ are pairwise non-adjacent in $G$.

    We distinguish two subcases depending on the number of leaves of $T$.

    \subparagraph{Subcase A.1: $T$ has at least $k+1$ leaves.}
    Let $L(T) \subseteq M^{\textrm{down}}$ be the set of leaves of $T$ and
    let $I(T)\subseteq V(G)$ be the set of all internal vertices of $T$, i.e., $I(T) = V(T) \setminus L(T)$. 
    Pick arbitrary $w^{\textrm{down}} \in L(T)$ and let $w^{\textrm{up}}$ be the neighbor of $w^{\textrm{down}}$ in $M$.
    Let $M'$ be the set of edges of $M$ that are incident to $L(T) \setminus \{w^{\textrm{down}}\}$. Note that $|M'| \geq k$.
    Finally, let $Z=I(T)\cup\{w^{\textrm{down}},w^{\textrm{up}}\}\cup \bigcup_{j<i-1} L_j$.
    We observe that $Z$ is connected and every vertex in $V(M')$ has a neighbor in $Z$:
    for vertices of $L_{i-1}$ such neighbor is in $L_{i-2}$ and for vertices in $L_i$ such neighbor is in $I(T)$.
    Thus, by \cref{obs:easywindmill}, $G$ contains the $k$-windmill as an induced minor, a contradiction.
    
    \subparagraph{Subcase A.2: $T$ has at most $k$ leaves.}
    We argue that in this case there is a path $P$ in $T$ with all internal vertices of degree 2 (in $T$) containing at least $2k$ vertices of $M^{\textrm{down}}$.
    To see why this holds, consider the tree $\widetilde{T}$ obtained from $T$ by contracting all vertices of degree two (in $T$).
    Since $\widetilde{T}$ has at most $k$ leaves and no degree-two vertices, it has at most $2k$ vertices and therefore at most $2k-1$ edges.
    Since $T$ spans at least $4k^2$ vertices of $M^{\textrm{down}}$, there is an edge of $\widetilde{T}$ corresponding to the claimed path $P$ in $T$.

    \begin{claim}    
        $P$ is an induced path in $G$.
    \end{claim}
    \begin{claimproof}
        Let the consecutive vertices of $P$ be $v_1,v_2,\ldots,v_p$ and suppose there are indices $j,j'$, where $j +1 < j'$, such that $v_jv_{j'}$ is an edge of $G$.
        If $v_{j+1},\ldots,v_{j'-1} \notin M^{\textrm{down}}$, then replacing $P$ in $T$ by $v_1,\ldots,v_j,v_{j'},\ldots,v_p$ we obtain a smaller tree that spans $M^{\textrm{down}}$ and has the same set of leaves as $T$, contradicting the minimality of $T$.
        So let $j'' \in [j+1,j''-1]$ be the smallest index for which $v_{j''} \in M^{\textrm{down}}$.
        Removing the edge $v_{j''}v_{j''+1}$ from $T$ and adding the edge $v_jv_{j'}$ instead, we obtain a tree with the same vertex set as $T$, but more leaves in $M^{\textrm{down}}$, contradicting the choice of $T$.        
    \end{claimproof}

    Since $P$ contains at least $2k$ vertices of $M^{\textrm{down}}$ and these vertices are pairwise nonadjacent, we can find a collection $\mathcal{Q}$ of $k$ pairwise vertex-disjoint and anticomplete subpaths of $P$, each having two vertices of $M^{\textrm{down}}$ as endpoints -- see~\cref{fig:caseA} for an illustration.
                
    Therefore, by removing all vertices of $G'$ that are not in a path from $\mathcal{Q}$, contracting each path in $\mathcal{Q}$ to a single edge, and contracting all edges in $\bigcup_{j< i}L_j$, we get the $k'$-windmill for $k' \geq k$ (and thus, also the $k$-windmill) as an induced minor, a contradiction.

    \newcommand{\tikzpa}[2]{
        \draw[smooth,#2] plot coordinates {(#1,0) (#1+0.25,-0.2) (#1+0.5,-0.1) (#1+0.7,-0.2) (#1+1,0)};
    }
    
    \begin{figure}[ht]
        \centering
        \begin{tikzpicture}[every node/.style={minimum size = 4pt, inner sep = 0pt}]
        
        \fill[gray!20!white] (0,0.5) -- (5,0.5) -- (2.5,2) -- (0,0.5);

        \node (L) at (2.5,1.2) {$\bigcup_{j< i}L_j$};
        
        \tikzpa{0}{black,line width=1pt}
        \tikzpa{1}{gray}
        \tikzpa{2}{black,line width=1pt}
        \tikzpa{3}{gray}
        \tikzpa{4}{black,line width=1pt}
        \foreach \x in {0,1,...,5} {
            \draw[blue] (\x,0) -- (\x,0.5);
            \node[circle,draw,fill=white] at (\x,0) {};
            \node[circle,draw,fill=white] at (\x,0.5) {};
            }

        \node (q) at (2.5,-1.1) {$\mathcal{Q}$}; 
        \draw[->] (2.2,-1) to[bend left=10] (0.8,-0.3);
        \draw[->] (2.8,-1) to[bend right=10] (4.2,-0.3);
        \draw[->] (2.5,-0.8) to (2.5,-0.3);
        
        \end{tikzpicture}
        \caption{An illustration of Case A.}
      \label{fig:caseA}
    \end{figure}
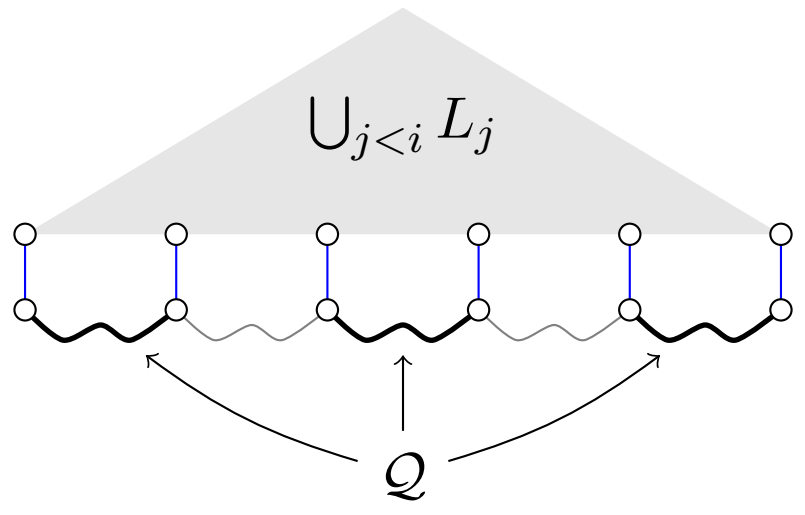
    This concludes Case A.

    \paragraph{Case B: Every component of $G'$ has fewer than $4k^2$ vertices of $M^{\textrm{down}}$.}
    Let $\cC$ be the set of components of $G'$ that intersect $M^{\textrm{down}}$.
    As $M$ has size at least $12k^3$, we have $|\cC|\geq 3k$.

    Our plan, towards finding a $k$-windmill as an induced minor and therefore contradicting our assumption for the existence of $M$, is the following.
    We will find a subset $\cC' \subseteq \cC$ of size $k$,
    and for every component $C \in \cC'$ we will find an edge $e_C=xy$ and two paths $P_x, P_y$ in $G$ such that:
    \begin{itemize}
        \item paths $P_x,P_y$ connect, respectively, $x$ and $y$ to the root $r$ and $y\notin V(P_x), x\notin V(P_y)$, and
        \item for $C' \in \cC' \setminus \{C\}$, the edge $e_{C'}$ is vertex-disjoint with both $P_x$ and $P_y$, 
        \item the set $\{e_C ~|~ C \in \cC'\}$ is an induced matching.
    \end{itemize}
    Note that this is sufficient to get the $k$-windmill as an induced minor.
    Indeed, we first remove from $G$ all vertices that do not belong to any of the chosen paths.
    Second, we contract all edges in the chosen paths, except for the edges incident to the vertices of the induced matching $\{e_C ~|~ C \in \cC'\}$.
    This way we obtain a $k$-windmill as an induced minor in $G$, a contradiction.

    It remains to define $\cC'$ and the special edge and paths for every component in $\cC'$.
    
    We distinguish three types of component in $\cC$:
    \begin{enumerate}[(I)]
        \item a single-vertex component,
        \item a tree with at least two vertices,
        \item not a tree.
    \end{enumerate}
    As $|\cC| \geq 3k$, there is a set $\cC' \subseteq \cC$ of size $k$ such that all components in $\cC'$ have the same type.
    We will consider three subcases, depending on the type of components in $\cC'$.

    \subparagraph{Subcase B.1: Type (I).}
    Let $v$ be the unique vertex in $C$, clearly $v \in M^{\textrm{down}}$.
    Let $u \in M^{\textrm{top}}$ be such that $vu \in M$.

    As $v$ has mininum degree 2 in $v$, there is $w \neq u$ such that $vw \in E(G)$.
    Note that $w \in L_{i-1}$.
    We pick $e_C = vu$.
    The path $P_u$ is a root path for $u$, 
    and the path $P_v$ is obtained by first following the edge $vw$ and then a root path of $P_w$. See~\cref{fig:caseBtypeI}.
    
    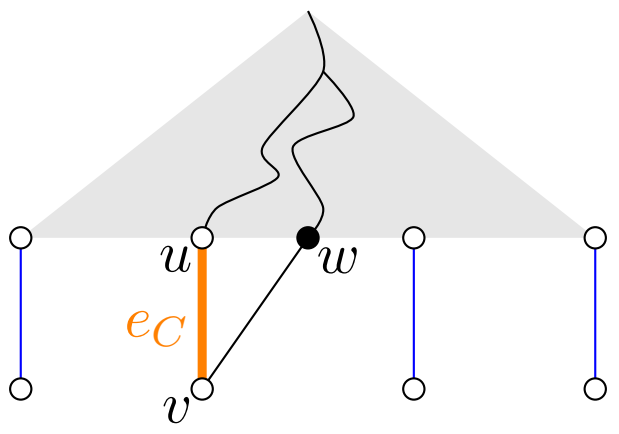
\begin{figure}[ht]
        \centering
        \begin{tikzpicture}[every node/.style={minimum size = 4pt, inner sep = 0pt}]
        
        \fill[gray!20!white] (0.6,0.5) -- (4.4,0.5) -- (2.5,2) -- (0.6,0.5);

        \draw[smooth] plot coordinates {(1.8,0.5) (1.9,0.7) (2.3,0.9) (2.2,1.1)  (2.6,1.6) (2.5,2)};  
        \node[circle,draw,fill=black,label={310:$w$}] at (2.5,0.5) {}; 
        \draw (1.8,-0.5) -- (2.5,0.5); 
        \draw[smooth] plot coordinates {(2.5,0.5) (2.6,0.7) (2.4,1.1) (2.8,1.3) (2.6,1.6)}; 
        
        \foreach \x in {0.6,3.2,4.4} {
            \draw[blue] (\x,-0.5) -- (\x,0.5);
            \node[circle,draw,fill=white] at (\x,-0.5) {};
            \node[circle,draw,fill=white] at (\x,0.5) {};
            }

        \draw[orange,line width=1.7pt] (1.8,-0.5) to (1.8,0.5);    
        \node[label={180:\textcolor{orange}{$e_C$}}] at (1.8,-0.1) {};\node[circle,draw,fill=white,label={220:$u$}] at (1.8,0.5) {};
        \node[circle,draw,fill=white,label={220:$v$}] at (1.8,-0.5) {};

        \end{tikzpicture}
        \caption{An edge $e_C$ for a component of type (I).}
      \label{fig:caseBtypeI}
    \end{figure}

    Now let us argue that such choice of $e_C$ and paths satisfies all required properties.
    The first property is obvious.
    The second and the third property follow from the fact that $\{e_C ~|~ C \in \cC'\} \subseteq M$ and $M$ is an induced $(L_{i-1},L_i)$-matching.

    \subparagraph{Subcase B.2: Type (II).}
    Consider any component $C \in \cC'$. Let $v$ be any vertex in $V(C) \cap  M^{\textrm{down}}$
    and $u$ be the unique neighbor of $v$ in $M^{\textrm{top}}$.
    Let $z$ be a leaf of $C$ that is distinct from $v$; it exists as $C$ is a tree on at least two vertices.
    Since $z$ has degree one in $C$ and degree at least two in $G$, it has a neighbor $z'$ in $L_{i-1}$.

    Let $e_C = vv'$ be the first edge on the unique $v$-$z$-path in $C$.
    We also set $P_v$ to be the path obtained by first following the edge $vu$ and then a root path of $u$,
    and $P_{v'}$ to be the path obtained by first following the $v'$-$z$-path in $C$, then the edge $zz'$, and then a root path of $z'$. See~\cref{fig:caseBtypeII}.

    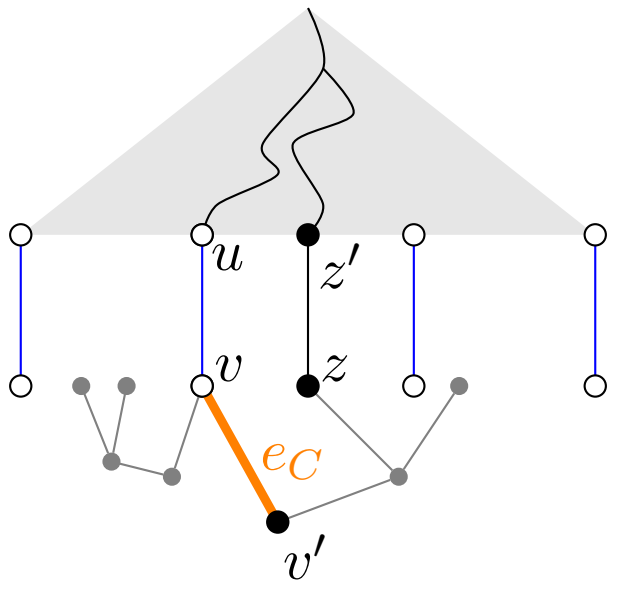
\begin{figure}[ht]
        \centering
        \begin{tikzpicture}[every node/.style={minimum size = 4pt, inner sep = 0pt}]
        
        \fill[gray!20!white] (0.6,0.5) -- (4.4,0.5) -- (2.5,2) -- (0.6,0.5);

        \draw[smooth] plot coordinates {(1.8,0.5) (1.9,0.7) (2.3,0.9) (2.2,1.1)  (2.6,1.6) (2.5,2)};  
        \node[circle,draw,fill=black,label={320:$z'$}] at (2.5,0.5) {}; 
        \draw (2.5,-0.5) -- (2.5,0.5); 
        \draw[smooth] plot coordinates {(2.5,0.5) (2.6,0.7) (2.4,1.1) (2.8,1.3) (2.6,1.6)}; 
        
        \foreach \x in {0.6,1.8,3.2,4.4} {
            \draw[blue] (\x,-0.5) -- (\x,0.5);
            \node[circle,draw,fill=white] at (\x,-0.5) {};
            \node[circle,draw,fill=white] at (\x,0.5) {};
            }

        \node[circle,draw=gray,fill=gray,minimum size = 3pt] (s) at (1.6,-1.1) {};
        \node[circle,draw=gray,fill=gray,minimum size = 3pt] (w) at (1.2,-1) {};
        \node[circle,draw=gray,fill=gray,minimum size = 3pt] at (1,-0.5) {};
        \node[circle,draw=gray,fill=gray,minimum size = 3pt] at (1.3,-0.5) {};
        \node[circle,draw=gray,fill=gray,minimum size = 3pt] (l) at (3.1,-1.1) {};
        \node[circle,draw=gray,fill=gray,minimum size = 3pt] (t) at (3.5,-0.5) {};
        \draw[gray]
            (s) to (1.8,-0.5)
            (s) to (1.2,-1)
            (w) to (1,-0.5) 
            (w) to (1.3,-0.5)
            (l) to (2.3,-1.4)
            (l) to (2.5,-0.5)
            (l) to (t)
            ;
        
        \draw[orange,line width=1.7pt] (1.8,-0.5) to (2.3,-1.4);
        \node[label={right:\textcolor{orange}{$e_C$}}] at (2.1,-1) {};
        
        \node[circle,draw,fill=white,label={310:$u$}] at (1.8,0.5) {};
        \node[circle,draw,fill=white,label={20:$v$}] at (1.8,-0.5) {};
        \node[circle,draw,fill=black,label={20:$z$}] at (2.5,-0.5) {};
        \node[circle,draw,fill=black,label={290:$v'$}] at (2.3,-1.4) {};

        \end{tikzpicture}
        \caption{An edge $e_C$ for a component of type (II).}
      \label{fig:caseBtypeII}
    \end{figure}

    As $e_C$ is contained in $C$, all required properties are clearly satisfied.

    \subparagraph{Subcase B.3: Type (III).}
    Consider any component $C \in \cC'$. Again, let $v$ be any vertex in $V(C) \cap  M^{\textrm{down}}$
    and $u$ be he unique neighbor of $v$ in $M^{\textrm{top}}$.

    Fix an induced cycle contained in $C$.
    Let $Q$ be a shortest path in $C$ joining $v$ with a vertex of the cycle; call this vertex $z$.
    We set $e_C$ to be any edge $xy$ of the cycle that does not contain $z$.
    The paths $P_x$ and $P_y$ are obtained by first taking two internally vertex-disjoint paths along the cycle connecting $x$ and $y$ to $z$, then following the path $Q$ to $v$, then the edge $vu$, and then a root path of $u$. See~\cref{fig:caseBtypeIII}.

        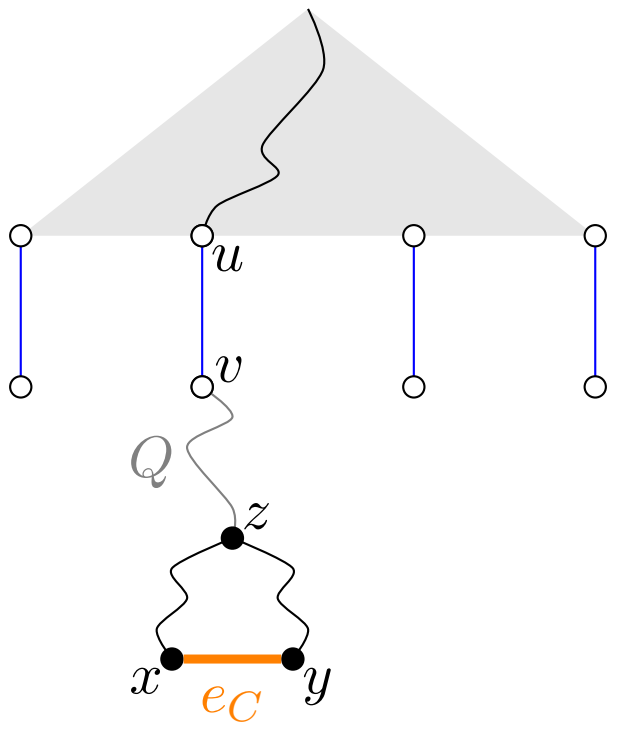
\begin{figure}[ht]
        \centering
        \begin{tikzpicture}[every node/.style={minimum size = 4pt, inner sep = 0pt}]
        
        \fill[gray!20!white] (0.6,0.5) -- (4.4,0.5) -- (2.5,2) -- (0.6,0.5);

        \draw[smooth] plot coordinates {(1.8,0.5) (1.9,0.7) (2.3,0.9) (2.2,1.1)  (2.6,1.6) (2.5,2)};  

        \foreach \x in {0.6,1.8,3.2,4.4} {
            \draw[blue] (\x,-0.5) -- (\x,0.5);
            \node[circle,draw,fill=white] at (\x,-0.5) {};
            \node[circle,draw,fill=white] at (\x,0.5) {};
            }

        \node[circle,draw,fill=black,label={south west:$x$}] (x) at (1.6,-2.3) {};
        \node[circle,draw,fill=black,label={south east:$y$}] (y) at (2.4,-2.3) {};
 
        \draw[smooth,gray] plot coordinates {(1.8,-0.5) (2,-0.7) (1.7,-0.9) (2,-1.3) (2,-1.5)}; 
        \node[label={left:\textcolor{gray}{$Q$}}] at (1.7,-1) {};

        \draw[orange,line width=1.7pt] (x) to (y);
        \node[label={below:\textcolor{orange}{$e_C$}}] at (2,-2.4) {};
        
        \draw[smooth] plot coordinates {(x) (1.5,-2.1) (1.7,-1.9) (1.6,-1.7) (2,-1.5)}; 
        \draw[smooth] plot coordinates {(y) (2.5,-2.1) (2.3,-1.9) (2.4,-1.7) (2,-1.5)}; 

        \node[circle,draw,fill=white,label={310:$u$}] at (1.8,0.5) {};
        \node[circle,draw,fill=white,label={20:$v$}] at (1.8,-0.5) {};
        \node[circle,draw,fill=black,label={north east:$z$}] at (2,-1.5) {};

        \end{tikzpicture}
        \caption{An edge $e_C$ for a component of type (III).}
      \label{fig:caseBtypeIII}
    \end{figure}
    
    Similarly as in the previous case, $e_C$ is contained in $C$, so all required properties are satisfied.

    \medskip
    This completes the case analysis and the proof of the claim.
\end{claimproof}

Now, the statement of the theorem follows from \cref{lem:from-bfs-to-tree}.
\end{proof}

\subsection{Dichotomy for induced minors}
Equipped with \cref{thm:windmill}, we can present the following dichotomy result.

\thmdichoindminordeg*
\begin{proof}
    The cases that define classes of bounded induced matching treewidth follow directly from combining \cref{thm:dichoindminor} and \cref{thm:windmill}.
    Thus, let us assume that $H$ is not one of these graphs.
    By \cref{thm:dichoindminor} and \cref{obs:afterindminor}, we may assume that $H$ contains a $k$-windmill as an induced subgraph for some $k \geq 2$.
    
    Fix an induced copy $W$ of a $k$-windmill with $k$ largest possible.
    Let $x$ be the central vertex, and all other vertices are referred to as \emph{satellites}.
    Each triangle of this windmill is called a \emph{sail}, i.e., each sail contains two adjacent satellites and $x$.
    
    Since $H$ is not a windmill, there is a vertex not in $W$.
    Let us consider two cases.

    \paragraph{Case A. There is a vertex $v$ non-adjacent to $x$.}
    If there is a sail of $W$ that is non-adjacent to $v$,
    then $H$ contains $K_3+P_1$ as an induced minor, and the conclusion follows from \cref{lem:min-unbounded}~\ref{it:join}.

    Thus, $v$ is adjacent to all sails of $W$.
    Let $H'$ be the subgraph of $H$ induced by the vertices of $W$ and $v$.
    Note that $H'$ is 2-connected and contains an induced $P_2 + P_1$, consisting of two satellites from one sail and one satellite from another sail.
    Thus, by \cref{lem:min-unbounded}~\ref{it:P2P1},
    induced matching treewidth is unbounded in $H$-induced-minor-free graphs, even restricted to graphs of minimum degree at least 2.

    \paragraph{Case B. Every vertex of $H$ is adjacent to $x$.}
    If $H-x$ is connected, then $H$ is 2-connected, and since it contains an induced $P_2 + P_1$, the conclusion follows from \cref{lem:min-unbounded}~\ref{it:P2P1}.
    Thus, $H-x$ has at least two components.
    Furthermore, as $H$ is not an induced minor of a windmill, there is a component of $H-x$ that has at least three vertices, i.e., it contains either $P_3$ or $K_3$.
    Summing up, $H$ contains $P_3+P_1$ or $K_3+P_1$ as an induced minor, and the conclusion follows from \cref{lem:min-unbounded}~\ref{it:join}.
    This completes the proof.
\end{proof}

Similarly to \cref{obs:afterindminor} we note that the only case where the graphs in our construction are not 2-connected, i.e., where we use \cref{lem:min-unbounded}~\ref{it:P2P1}, is the case that $H$ contains a 2-connected induced subgraph that contains $P_2 + P_1$ as an induced subgraph.
Thus, we obtain the following refinement of \cref{obs:afterindminor}.

\begin{observation}\label[observation]{obs:aftermindeg}
Let $H$ be a graph for which $H$-induced-minor-free graphs have unbounded induced matching treewidth.
Induced matching treewidth in $H$-induced-minor-free graphs remains unbounded, even restricted to 2-connected graphs, unless $H$ contains a 2-connected induced subgraph that contains $P_2 + P_1$ as an induced subgraph.
\end{observation}

\newpage
\section{Restricting to 2-connected graphs}\label{sec:2connected}
In this section we restrict the considered graphs even more, by assuming that they are 2-connected.
Recall that this assumption does not change anything for tree-independence number (\cref{lem:block-reduction}), nor for induced matching treewidth of graphs that exclude a fixed graph $H$ as a subgraph, induced subgraph, a minor, a topological minor, or an induced subdivision.

Again it appears that the remaining case, i.e., $H$-induced-minor-free graphs, behaves differently.

\subsection{New bounded case}

\begin{figure}[ht]
    \centering
    \begin{tikzpicture}[every node/.style={draw, circle, fill = black, minimum size = 2.5pt, inner sep = 0pt}]

        \node (A) at (0,1) {};
        \node (B) at (1,1) {};
        \node (C) at (1,0) {};
        \node (D) at (0,0) {};
        \node (H) at (0.5,1.5) {};

        \draw (A) -- (B) (B) -- (C) (C) -- (D) (D) -- (A) (A) -- (C) (B) -- (D);
        \draw (A) -- (H) -- (B);
    \end{tikzpicture}
    \caption{The graph $\widehat{K}_4$.}
    \label{fig:hatk4}
\end{figure}
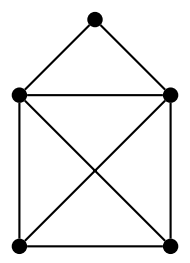
First, we exhibit one new graph whose exclusion bounds induced matching treewidth.

\begin{theorem}\label{thm:2con-k4hat}    
    Let $G$ be a $2$-connected graph that excludes $\widehat{K}_4$ as an induced minor. Then $\yolov(G) \leq 4$.
\end{theorem}

To prove~\Cref{thm:2con-k4hat}, we use the following result of B{\l}asiok, Kami\'nski, Raymond, and Trunck~\cite[Theorem 2]{BlasiokKTT19}.
Here, a \emph{prism} is the complement of $C_6$, i.e., two disjoint triangles joined with a perfect matching.

\begin{theorem}[B{\l}asiok, Kami\'nski, Raymond, and Trunck~\cite{BlasiokKTT19}]\label{thm:k4hat-structure}
    Let $G$ be a $2$-connected graph that excludes $\widehat{K}_4$ as an induced minor. Then one of the following holds:
    \begin{enumerate}
        \item $G$ is $K_4$-minor-free;
        \item $G$ is a subdivision of $K_4$, of $K_{3,3}$, or of the prism;
        \item $V(G)$ can be partitioned to sets $W,M\subseteq V(G)$ such that $G[W]$ is a wheel on at most five vertices and $G[M]$ is a complete multipartite graph;
        \item $V(G)$ can be partitioned to sets $C,I\subseteq V(G)$ such that $G[C]$ is a cycle, $G[I]$ is edgeless, and every vertex of $I$ has the same neighborhood on $C$.
    \end{enumerate}
\end{theorem}

Now let us proceed to the proof of ~\Cref{thm:2con-k4hat}
\begin{proof}[Proof of~\Cref{thm:2con-k4hat}]
We examine all cases given by Lemma~\ref{thm:k4hat-structure}.
\begin{enumerate}[label={Case \arabic*.}]
    \item If $G$ is $K_4$-minor-free, then $\tw(G) \leq 2$~\cite{Duffin:1965:TopologyOfSeriesParallelNetworks}, and thus, $\yolov(G) \leq 3$. 
    \item It is easy to verify that the treewidth of $K_4$, $K_{3,3}$, or a prism is 3.
    Since subdividing edges does not increase treewidth, we conclude that $\tw(G) \leq 3$ and thus, $\yolov(G) \leq 4$.
    \item Let $F$ be an maximum-size induced matching in $G$.
    At most one edge of $F$ might be contained in $W$ and at most one edge might be contained in $M$.
    Moreover, there might be at most 2 edges of $F$ with one endvertex in $M$ and the other in $W$, as $W$ does not contain a tree-element independent set. Consequently, $|F|\le 4$, and thus, $\yolov(G) \leq 4$.
    \item Let $\cT$ be a tree decomposition of $G[C]$ with width 2, and modify it by adding every vertex from $I$ to every bag.
    It is straightforward to verify that this way we obtain a tree decomposition of $G$.
    Consider any bag $X_u$ of this decomposition. Note that $|X_u \cap C| \leq 3$ and since $I$ is independent,
    the size of a largest induced matching touching $X_u$ is at most 3.
    Consequently, $\yolov(G) \leq 3$.    
\end{enumerate}
This completes the proof.
\end{proof}

\subsection{Towards the dichotomy for 2-connected graphs}

While we are not able to provide a full dichotomy in this case, we show the following \emph{almost dichotomy} theorem.

\thmindminortwoconn*
\begin{proof}
    The bounded cases follow from the combination of \cref{thm:dichoindminordeg} and \cref{thm:2con-k4hat}.
    So let us consider a graph $H$ that is not one of these cases, nor  $\widehat{K}_4^{(q)}$ for some fixed $q$.
    
    By \cref{obs:aftermindeg} we may assume that some 2-connected induced subgraph of $H$ contains $P_1+P_2$ as an induced subgraph.
    Fix one copy of such subgraph $H'$ and an induced copy of $P_1+P_2$ in $H'$.
    Let $a,b,c$ be the vertices of this copy, where $bc \in E(H)$.

    As $H$ is not a subgraph of a matching, there are vertices in $V(H) \setminus \{a,b,c\}$.
    Note that if there is $v \in V(G) \setminus \{a,b,c\}$ that is adjacent to $\{b,c\}$ but non-adjacent to $a$, then $H$ contains either $K_3+P_1$ or $P_3+P_1$ and we are done by \cref{lem:min-unbounded}~\ref{it:join}. Thus, every vertex in $V(H) \setminus \{a,b,c\}$ that is adjacent to $\{b,c\}$ is also adjacent to $a$.

    Next, if there is $v$ adjacent to, say $b$ (and thus to $a$), but not $c$, then $H$ contains an induced $P_4$ and again we are done by \cref{lem:min-unbounded}~\ref{it:join}.

    Since $H'$ is 2-connected, there are two internally disjoint paths from $a$ to $\{b,c\}$ in $H'$ and thus in $H$.
    Pick one such path, without loss of generality we can assume that it is induced.
    If the path has more than one internal vertex, then $H$ contains an induced $P_4$ and we are done. So each of the two paths has exactly one internal vertex.
    Thus, there are two vertices $x,y$, both adjacent to all $\{a,b,c\}$.

    If $xy \notin E(H)$, then $H[\{a,b,c,x,y\}]$ is isomorphic to $K'_4$, the graph obtained from $K_4$ by once subdividing an edge. Unboundedness of induced matching treewidth in $H$-induced-minor-free graphs follows from \cref{lem:min-unbounded}~\ref{it:k4}.
    
    So, let us assume that $xy \in E(H)$.
    This means that the set $\{a,b,c,x,y\}$ induces a copy of  $\widehat{K}_4$ in $H$.

    If there is a vertex in $V(H) \setminus \{a,b,c,x,y\}$ adjacent to $a$,
    then $a,b,c$ are of degree at least 3 and they induce $P_2+P_1$ in $H$, so we are done by \cref{lem:min-unbounded}~\ref{it:P2P1-mindeg3}.

    So suppose that every vertex  $v \in V(H) \setminus \{a,b,c,x,y\}$ is non-adjacent to $a$.
    If there is $v\in V(H) \setminus \{a,b,c,x,y\}$ that is adjacent to at least one of $\{b,c\}$, then $H[\{v,b,c,a\}$ is either $K_3 + K_1$ or $P_3 + K_1$ -- in both cases we are done by  by \cref{lem:min-unbounded}~\ref{it:join}.

    Suppose there is a vertex $v$ whose neighborhood in $\{b,c,a,x,y\}$ is not $\{x,y\}$.
    As $v$ is non-adjacent to $a,b,c$, by symmetry let us assume that either $v$ has no neighbors in $\{b,c,a,x,y\}$, or its only neighbor in this set is $x$.
    But then, $H[\{v,b,c,y\}]$ is isomorphic to $K_3 + K_1$ and we are done by  by \cref{lem:min-unbounded}~\ref{it:join}.

    Consequently, every vertex $v$ in $V(H) \setminus \{a,b,c,x,y\}$ is adjacent to both $\{x,y\}$ and none of $\{b,c,a\}$.
    Finally, suppose that  $V(H) \setminus \{a,b,c,x,y\}$ is not independent, i.e., it contains two adjacent vertices $v,v'$.
    Then $H[\{v,v',x,a\}$ is $K_3 + K_1$ and we are done by  by \cref{lem:min-unbounded}~\ref{it:join}.

    This means that $H$ is isomorphic to $\widehat{K}_4^{(q)}$ for some $q$. This contradiction completes the proof.
    \end{proof}

\newpage
\section{Algorithmic consequences}\label{sec:algos}
Boundedness of induced matching treewidth in a certain class of graphs allows us to solve numerous computational problems in polynomial time.
Furthermore, for many problems we can reduce solving them to blocks of the input graph.
Thus, we can make use of the positive statement in \cref{thm:indminortwoconn}.

\subsection{Finding large induced sparse subgraphs.}

For a fixed integer $w$ and fixed \textsf{CMSO}$_2$ sentence\footnote{A \textsf{CMSO}$_2$ sentence is a formula in the second-order logic of graphs with quantifiers over vertices and edges and without free variables. For more details, see, e.g.,~\cite{platypus}.} $\psi$, by $(\tw \leq w,\psi)$-\MWIS we denote the problem where we are given a graph $G$ with vertex-weights, and we want to find a maximum-weight subset $X \subseteq V(G)$ such that the graph induced by $X$ has treewidth at most $w$ and satisfies the sentence $\psi$, or report that no such set exists.

This formalism captures many classic computational problems, like \textsc{Maximum Weight Independent Set} (for $w=0$ and $\psi$ being a sentence satisfied by all graphs),  \textsc{Feedback Vertex Set} (by complementation, for $w=1$ and $\psi$ being a sentence satisfied by all graphs).

By  $(\tw \leq w,\psi)$-\MWISa we denote the generalization of  $(\tw \leq w,\psi)$-\MWIS  where vertices of the input graph might be colored using constant-size set of colors the formula $\psi$ 
uses these colors.
The following result was shown by Bodlaender, Fomin, Korhonen~\cite{BodlaenderFK26} (see the discussion in the conclusion).

\begin{theorem}[Bodlaender, Fomin, Korhonen~\cite{BodlaenderFK26}]\label{thm:cmsomwis}
    For every $w$ and $\psi$, the \emph{$(\tw \leq w,\psi)$-\MWISa} problem can be solved in polynomial time on graphs of bounded induced matching treewidth.
\end{theorem}

Next, we show that in order to solve $(\tw \leq w,\psi)$-\MWISa in some hereditary class of graphs,
it is sufficient to focus on 2-connected graphs from this class. The proof is a standard dynamic-programming process along block-cut trees.

\begin{lemma}\label[lemma]{lem:cmsomwis-2conn}
    Let $\cG$ be a hereditary class of graphs (whose vertices are colored with a fixed number of colors) and let $\cG_{\textrm{2-conn}}$ be the family of 2-connected members of $\cG$.
    If, for every fixed $r$ and $\psi$,
    and every graph in $\cG_{\textrm{2-conn}}$,
    the  \emph{$(\tw \leq w,\psi)$-\MWISa} problem can be solved in polynomial time, then it can also be solved in polynomial time for every graph in $\cG$.
\end{lemma}
\begin{proof}
    Let $\ell$ denote the (fixed) number of colors appearing in graphs from $\mathcal{C}$.
    The algorithm will treat graphs with one distinguished vertex, which we will mark by an extra fresh color.
    Let $q$ be the \emph{quantifier rank} of $\psi$, i.e., the maximum number of nested quantifiers in it.
    It is well-known that there are only finitely many syntactically non-equivalent \textsf{CMSO}$_2$ sentences over graphs (colored with $\ell+1$ colors) and of quantifier rank at most $q$. We denote this set of formulas by $\Phi_{q,\ell+1}$ and note that there is a sentence $\psi'$ in $\Phi_{q,\ell+1}$ that is equivalent to $\psi$, in the sense that for every graph $G$, $G$ satisfies $\psi$ if and only if it satisfies $\psi'$. Since $\psi$ is fixed, we may assume that $\Phi_{q,\ell}$ and $\psi'$ are also available to our algorithm.

    Given a graph $G\in\cG$ colored with $\ell$ colors, we assume that an arbitrary vertex is marked with an extra ($(\ell+1)$-th) color. Our goal is to prove the statement of the lemma for every sentence in $\Phi_{q,\ell+1}$, i.e., for each $\phi\in\Phi_{q,\ell+1}$, find a maximum-weight subset $X\subseteq V(G)$ such that the graph induced by $X$ has treewidth at most $w$ and satisfies $\phi$, or report that no such set exists.
    We call such a set a \emph{maximum-weight solution} for $\phi$ in $G$.

    We first show how to deal with the case where $G$ is connected and then we lift the argument to disconnected graphs. If $G$ is 2-connected, we are done by the assumption.
    Thus, suppose otherwise and let $T$ be the block-cut tree of $G$.
    For every non-root block $b$ of $T$, we use $c_b$ to denote the cutvertex between $b$ and its parent block. 
    We proceed by dynamic programming, using the structure of $T$ in a bottom-up fashion. 

    Let us first establish some notation. For a node $u$ of $T$, we denote by $G_u$ the subgraph of $G$ induced by the blocks/cutvertices in the subtree of $T$ rooted at $u$. The distinguished vertex (colored by the extra color) of $G_u$ is $u$ if $u$ is a cutvertex, or the cutvertex $c_b$, in case $u$ is a block (in case $b$ is the root block, we set $c_b$ be an arbitrarily chosen vertex in $b$).    
    
    For each node $u$ of the block-cut tree and each formula $\phi\in\Phi_{q,\ell}$, we will store an entry $D_{\phi}(u)$ consisting of a number $p$, corresponding to the weight of a maximum-weight solution $S$ for $\phi$ in $G_u$, as well as the part of the set $S$ in the block/cutvertex $u$ and pointers to entries of the form $D_{\phi'}(u')$, where $\phi'\in\Phi_{q,\ell}$ and $u'$ is a child of $u$. The latter information will be used so as to recursively recover the maximum-weight solution for $G_u$. If a solution does not exist for certain $\phi$, we set $D_{\phi}(u)=\bot$.
    

    We start from leaf blocks. Note that a leaf block $b$ either consists of two vertices or it is 2-connected. Therefore, for each $\phi\in\Phi_{q,\ell}$, we find a maximum-weight solution for $\phi$ in $G_b$ (if such exists) by either using brute-force (in the former case), or calling the algorithm for 2-connected instances (in the latter case). In both cases, we store the corresponding information in $D_{\phi}(b)$.
    
    For a non-leaf block $b$ for which we have already $D_{\phi'}(b')$ for every $\phi'\in\Phi_{q,\ell}$ and every child block $b'$ of $b$, we proceed as follows.
    We use the following fact about \textsf{CMSO}$_2$, following from the well-known ``Feferman-Vaught'' compositionality of types of \textsf{CMSO}$_2$ (see~\cite{CourcelleEngelfriet2012}).\medskip

    \emph{Fact:} Let $\phi\in\Phi_{q,\ell+1}$. There is a collection of pairs of sentences $(\psi_1^j,\psi_2^j)_{j\in\{1,\ldots,p\}}\subseteq \Phi_{q,\ell+1}\times\Phi_{q,\ell+1}$ for a $p\in\mathbb{N}$ depending only on $q$ and $\ell$, such that the following holds.
    Let $G_1,G_2$ be two graphs (colored with $\ell$ colors) with a distinguished vertex $v_1$ and $v_2$ in each of them respectively (colored with an ($\ell$+1)-th color). Also let $\hat{G}$ be the graph obtained by taking the disjoint union of $G_1$ and $G_2$ and identifying $v_1$ with $v_2$. Then for every set $S\subseteq V(\hat{G})$,
    $\hat{G}[S]$ satisfies $\phi$ if and only if there is a $j\in\{1,\ldots,p\}$ and sets $S_1\subseteq V(G_1)$ and $S_2\subseteq V(G_2)$ such that for each $i\in\{1,2\}$, $G_i[S_i]$ satisfies $\psi_i^j$ and $S=S_1\cup S_2$. We call such a pair of formulas \emph{admissible for $\phi$}.\medskip

    In the above fact, it is easy to observe that the treewidth of $\hat{G}[S]$ is the maximum of the treewidth of $G_1[S_1]$ and $G_2[S_2]$ (since these two graphs have at most one vertex in common) and that every maximal-weight set $S$ as above is the union of two maximal-weight sets $S_1$ and $S_2$ as above and vice versa.

    We now apply the following intermediate step: we combine solutions of children blocks with a common cutvertex.
    Let $c$ be a cutvertex of $G$ that belongs to $b$ and is not $c_b$ and let $b_1,\ldots,b_m$ be the children blocks of $b$ that share the cutvertex $c$ with $b$.
    We can compute $D_{\phi}(c)$ for every sentence $\phi\in\Phi_{q,\ell}$, by iteratively applying the above fact and computing (in time $\mathcal{O}(m)$) a collection of $m$-tuples of admissible formulas $(\phi_1^j,\ldots,\phi_m^j)_{j\in\{1,\ldots,z\}},$ where $z$ is a constant depending only on $\ell$ and $q$, such that $D_{\phi}(c)$ is determined by $D_{\phi_i^j}(b_i)$ for $i\in\{1,\ldots,m-1\}$, for some choice of $j\in\{1,\ldots,z\}$ maximizing the total weight of a solution.
    
    We next proceed to compute $D_{\phi}(b)$ for every $\phi\in\Phi_{q,\ell+1}$.
    Let $c_1,\ldots,c_m$ be the cutvertices of $G$ belonging to $b$ that are not $c_b$.
    In the subgraph $B$ of $G$ induced by the vertices in the block $b$, we assign a fresh color to each cutvertex $c_i$, corresponding to the formulas $\phi'\in\Phi_{q,\ell}$ for which $D_{\phi'}(c_i)\neq\bot$.
    Note that the number of colors in this new colored (2-connected) graph, denoted by $\tilde{B}$, is some constant $\ell'$ depending only on $q$ and $\ell$. 
    Also, the iterative application of the above fact implies that there is
    \begin{itemize}
        \item a  sentence $\tilde{\phi}$ with quantifier rank at most $q$, on the vocabulary of $\ell'$ colored graphs (which can be computed from $\phi$ in constant time) and
        \item a collection of $m$-tuples of admissible formulas $(\tau_1^j,\ldots,\tau_m^j)_{j\in\{1,\ldots,z'\}},$ where $z'$ is a constant depending only on $\ell$ and $q$ (which can be computed in $\mathcal{O}(m)$ time), 
    \end{itemize}
    such that a maximum-weight solution for $\phi$ in $G_b$ is obtained by the maximum-weight solution for $\tilde{\phi}$ in $\tilde{B}$ expanded by the recovered solution from $D_{\tau_1^j}(c_1),\ldots,D_{\tau_m^j}(c_m)$ for the tuple $(\tau_1^j,\ldots,\tau_m^j)$ that maximizes the total weight. This concludes the description of recursive steps along the block-cut tree.

    In case the input graph $G$ is disconnected, we first solve the problem in each connected component by applying the above procedure to compute maximum-weight solutions for every formula in $\Phi_{q,\ell+1}$ in every connected component. Then, by the (standard) analogue of the above fact applied for disjoint unions, this gives a maximum-weight solution in $G$ for every formula in $\Phi_{q,\ell+1}$.
    

    The overall running time is polynomial, provided that the algorithm for 2-connected instances runs in polynomial time.
\end{proof}

Combining \cref{thm:cmsomwis}, \cref{lem:cmsomwis-2conn}, and \cref{thm:indminortwoconn}, we obtain the following.
    \begin{restatable}{corollary}{corcmsomwis}
        For every $r$ and $\psi$, the \emph{$(\tw \leq r,\psi)$-\MWISa} problem can be solved in polynomial time when restricted to graphs from the following classes:
    \begin{multicols}{2}
    \begin{itemize}
        \item $K_{2,q}$-induced-minor-free for any fixed $q$,
        \item $W_4$-induced-minor-free,
        \item $K_5^-$-induced-minor-free,
        \item $q$-windmill-induced-minor-free for any fixed $q$,
        \item $\widehat{K}_4$-induced-minor-free.
    \end{itemize}
    \end{multicols}
\end{restatable}

\subsection{Finding large induced \texorpdfstring{$H$}{H}-colorable subgraphs.}
A \emph{homomorphism} from a graph $G$ to a graph $H$ is an edge-preserving mapping from $V(G)$ to $V(H)$.
For a fixed simple graph $H$, in the \textsc{Max Induced $H$-Colorable Subgraph} problem,
the input consists of a triple $(G,L,\rev)$,
where $G$ is a graph,
$L: V(G) \to 2^{V(H)}$ is a \emph{list function},
 and $\rev : V(G) \times V(H) \to \mathbb{Q}$ is a \emph{revenue function}.
A feasible solution for such an instance is a pair $(X,\phi)$, where $X \subseteq V(G)$ and $\phi$ is a homomorphism from $G[X]$ to $H$, respecting lists $L$.
The revenue of such a solution is given by $\sum_{v \in X} \rev(v,\phi(v))$.
Our goal is to find a solution of maximum possible revenue.

This formalism captures many classic problems line \MWIS, \textsc{Odd Cycle Transversal}, or \textsc{(List) $k$-Coloring}. 
We remark that the lists $L$ in the definition of the problem are somewhat redundant, as they can be modeled by adjusting the revenue function. However, we keep them as they make the description of the algorithm simpler.
We emphasize that the list of a vertex might be empty, then it cannot be included in any feasible solution at all.

The work of Yolov~\cite{Yolov18} yields the following result.

\begin{theorem}[Yolov~\cite{Yolov18}]\label{thm:Hcolorable}
    For every fixed simple graph $H$ and every fixed integer $k$,
    \textsc{Max Induced $H$-Colorable Subgraph}
    can be solved in polynomial time when resricted to graphs of induced matching treewidth at most $k$.
\end{theorem}

The following lemma shows that, while solving \textsc{Max Induced $H$-Colorable Subgraph}, it is sufficient to focus on 2-connected graphs. 

\begin{lemma}\label[lemma]{lem:Hcolorable-2conn}
    Let $\cG$ be a hereditary class of graphs such that and let $\cG_{\textrm{2-conn}}$ be the family of 2-connected members of $\cG$.
    If, for every graph in $\cG_{\textrm{2-conn}}$, the \textsc{Max Induced $H$-Colorable Subgraph} problem can be solved in polynomial time, then it can also be solved in polynomial time for every graph in $\cG$.
\end{lemma}
\begin{proof}
    Consider an instance $(G,L,\rev)$ of \textsc{Max Induced $H$-Colorable Subgraph}, where $G \in \cG$.
    Note that without loss of generality we may assume that $G$ is connected.
    If $G$ is 2-connected, we are done by the assumption.
    Thus, suppose otherwise and let $T$ be the block-cut tree of $G$.
    Let $r$ be an arbitrary leaf; recall that $r$ is a block.
    We treat $T$ as rooted at $r$; this defines parent-child relation between adjacent nodes.
    Let $v_r$ be any vertex of $r$ that does not belong to any other block (i.e., it is not the unique cutvertex that belongs to $r$).
    
    For a block node $b \neq r$, let $v_b$ be the cutvertex that is the parent of $b$ in $T$,
    and  the \emph{parent block} is the parent of $v_b$ in $T$.
    Analogously we define \emph{children blocks}.

    We proceed by dynamic programming, using the structure of $T$ in a bottom-up fashion.
    For each block $b$, and each $c \in L(v_b) \cup \{\bot\}$, we will compute the maximum possible revenue $\textsf{Opt}(b,c)$ of a solution $(X,\phi)$ in the subgraph induced by the blocks in the subtree of $T$ rooted at $b$ such that
    \[
    \begin{cases}
        v_b \in X \text{ and } \phi(v_b) = c & \text{ if } c \in L(v_b)\\
        v_b \notin X & \text{ if } c = \bot.
    \end{cases}
    \]
    Clearly, the optimum revenue is $\max_{c \in L(v_r) \cup \{\bot\}} \textsf{Opt}(r,c)$.
        
    Thus, fix some block $b$ and some $c \in L(v_b) \cup \{\bot\}$.
    Suppose first that $b$ is a leaf block.
    If $b$ has two vertices, then  $\textsf{Opt}(b,c)$ can be computed by brute force.
    Otherwise, the $b$ is 2-connected and $\textsf{Opt}(b,c)$ can be computed  by calling the algorithm for 2-connected instances
    with the list of $v_b$ set to $\{c\}$ if $c \in L(v_b)$ or to $\emptyset$ is $c = \bot$.

    Thus, suppose that $b$ is not a leaf block and $\textsf{Opt}(b',c')$ is computed for every child block $b'$ of $b$ and every $c' \in L(v_{b'}) \cup \{\bot\}$.
    First, we restrict the list of $v_b$ to $c$ or set it to $\emptyset$, as before.
    Next, for every child block $b'$ of $b$ and every $c' \in L(v_{b'})$,
    we set the revenue of the pair $(v_{b'},c')$ to $\textsf{Opt}(b',c') - \textsf{Opt}(b',\bot)$.
    It is straightforward to verify that $\textsf{Opt}(b,c)$ is equal to the revenue of the optimum solution to this instance plus $\sum_{b' \text{: child of } b} \textsf{Opt}(b',\bot)$.
    The overall running time is polynomial, provided that the algorithm for 2-connected instances runs in polynomial time.
\end{proof}

Combining \cref{thm:Hcolorable}, \cref{lem:Hcolorable-2conn}, and \cref{thm:indminortwoconn}, we obtain the following.
    \begin{restatable}{corollary}{corHcolorable}
    For every fixed graph $H$, the \textsc{Max Induced $H$-Colorable Subgraph} problem 
    can be solved in polynomial time when restricted to graphs from the following classes:
    \begin{multicols}{2}
    \begin{itemize}
        \item $K_{2,q}$-induced-minor-free for any fixed $q$,
        \item $W_4$-induced-minor-free,
        \item $K_5^-$-induced-minor-free,
        \item $q$-windmill-induced-minor-free for any fixed $q$,
        \item $\widehat{K}_4$-induced-minor-free.
    \end{itemize}
    \end{multicols}
\end{restatable}

\newpage
\section{Dominated balanced separators}\label{sec:dbs}
Let $\cG$ be a graph class.
We say that $\cG$ has \emph{$k$-dominated balanced separators} if for every $G \in \cG$ and every weight function $\wei : V(G) \to \RR_{\geq 0}$, there is $S \subseteq V(G)$ such that $|S| \leq k$ and every component of $G - N[S]$ has weight at most $\frac{1}{2} \wei(V(G))$.

First, we show that graphs of bounded induced matching treewidth admit $\Oh(1)$-dominated balanced separators.

\begin{lemma}\label[lemma]{lem:separators}
    The class of graphs of induced matching treewidth at most $t$ has $2t$-dominated balanced separators.    
\end{lemma}
\begin{proof}
    Let $G$ be a graph from $\cG$, $\wei : V(G) \to \RR_{\geq 0}$ be any weight function, and let $\cT = (T, \{X_u\}_{u \in V(T)})$ be its tree decomposition with $\mu(\cT,G) \leq t$.
    We aim to find a set $S$ with at most $2t$ vertices such that every component of $G - N[S]$ has weight at most $\frac{1}{2} \wei(V(G))$.
    
    From the properties of tree decompositions it follows that there exists $u \in V(T)$ such that every component of $G - X_u$ has weight at most $\frac{1}{2} \wei(V(G))$.
    Let $M$ be a largest induced matching in $G$ touching $X_u$; clearly $|M| \leq t$ and thus, $|V(M)| \leq 2t$.
    Notice that in $G - N[V(M)]$ there is no edge with at least one endpoint in $X_u$, as the existence of such an edge contradicts the maximality of $M$.    
    
    Thus, in $G - N[V(M)]$, there are components of two types: (i) components of $G - X_u$, and (ii) isolated vertices, all contained in $X_t$.
    Each component of type (i) has weight at most $\frac{1}{2} \wei(V(G))$ by the choice of $X_u$.
    If every component of type (ii) has weight at most $\frac{1}{2} \wei(V(G))$, we can set $S = V(M)$.
    If this is not the case, there is a unique vertex $v \in X_t \setminus N[V(M)]$ such that $\wei(v) > \frac{1}{2} \wei(V(G))$.
    Now it is sufficient to choose $S = \{v\}$.    
\end{proof}

Next, we show that in order to obtain $\Oh(1)$-dominated balanced separators in a hereditary class of graphs,
it is actually sufficient that its 2-connected members admit such separators.

\begin{lemma}\label[lemma]{lem:separators-2conn}
    Let $\cG$ be a hereditary class of graphs such that, for every $G \in \cG$, if $G$ is 2-connected, then it admits $t$-dominated balanced separators.
    Then $\cG$ has $t$-dominated balanced separators.
\end{lemma}
\begin{proof}    
    Let $G \in \cG$ and let $\wei : V(G) \to \RR_{\geq 0}$ be any weight function.
    Note that we can also assume that $G$ is connected -- otherwise we can focus on finding the balanced separator in the component of largest weight.
    
    Let $\cB$ be the set of blocks of $G$ and $\cC$ be the set of cutvertices.
    Consider the block-cut-tree $T_B$ of $G$: the vertex set of $T_B$ is $\cB \cup \cC$ and $vb$ is an edge if $b \in \cB$ and $v \in \cC$ and $v$ belongs to $b$.

    Define the weight $\wei_B : V(T_B) \to \RR_{\geq 0}$ as follows: for $v \in \cC$, we set $\wei_B(v) = \wei(v)$,
    and for $b \in \cB$ we set $\wei_B(b) = \sum_{v \in b \setminus \cC} \wei(v)$. Note that $\wei(V(G)) = \wei_B(V(T_B))$.
    
    Since $T_B$ is a tree, it has a node such that, once we remove it,  every component has weight at most $\frac{1}{2} \wei_B(V(T_B))$.
    If the separating node is a cutvertex or an edge, we pick these at most two vertices as a balanced separator in $G$.
    So assume that the separating node is a non-edge block, i.e., a maximal 2-connected subgraph of $G$, call it $G'$. 
    Clearly, $G'$ is an induced subgraph of $G$, so in particular it belongs to $\cG$.

    Consider components of $G - V(G')$.
    Each of them is adjacent to exactly one vertex of $G'$ (these vertices are cutvertices of $G$).
    For a vertex $v \in V(G')$, let $V_v$ denote the union of the vertex sets of all components of $G - V(G')$ that are adjacent to $v$.
    Clearly, if $v \notin \cC$, then $V_c = \emptyset$.

    Define a new weight function $\wei'$ on $G'$ as follows:
    for each $v \in V(G')$, we set $\wei'(v) = \wei(v) + \sum_{u \in V_v} \wei(u)$.
    Note that $\wei(V(G)) = \wei'(V(G'))$.
    
    By the assumption, there is $S \subseteq V(G')$ such that $|S| \leq t$ and every component of $G' - N_{G'}[S]$ has weight (according to $\wei'$) at most $\frac{1}{2} \wei'(V(G'))$.
    
    We claim that $N[S]$ is a balanced separator in $G$ with respect to $\wei$. Consider a component $C$ of $G - N[S]$.
    If $C$ intersects $V(G')$, then $\wei(V(C)) = \wei'(V(C) \cap V(G')) \leq \frac{1}{2} \wei(V(G')) = \frac{1}{2} \wei(V(G))$ as the weight of vertices of $V(G) \setminus V(G')$ was transferred to the vertices in $V(G') \cap \cC$.
    
    If $C$ does not intersect $V(G')$, then it is contained in one of the components of $G - V(G')$ and, by the choice of $G'$,
    we have $\wei(V(C)) \leq \frac{1}{2} \wei(V(G))$.

    This completes the proof.
\end{proof}

Combining \cref{thm:indminortwoconn} with \cref{lem:separators,lem:separators-2conn} yields the following corollary.

\cordbs*

\newpage
\section{Conclusion}\label{sec:outro}
Let us summarize the paper with pointing out some directions for further research.
First, we recall \cref{conj:fullhouses}.

\conjfullhouses*

Furthermore, we believe that the approach of preprocessing the graph and then asking for boundedness of some parameter is promising and not much explored~\cite{GaspersHP19,DabrowskiMNPR23,DBLP:journals/corr/abs-2004-05022}.
Possible preprocesing steps might include decomposing on clique cutsets~\cite{tarjan1985decomposition}, applying modular decomposition~\cite{habib2010survey} or contracting degree-2 vertices.

Finally, let us mention some further generalization of induced matching treewidth, introduced recently by Nikabadi and Rzążewski~\cite{amirwidth}, who were in turn inspired by a question of Bodlaender, Fomin, and Korhonen~\cite{BodlaenderFK26}.
For a family $\cH$ of graphs, the \emph{induced $\cH$-packing treewidth} is a tree-decomposition-based parameter where, for each bag of the decomposition, we measure the number of disjoint and pairwise nonadjacent copies of graphs from $\cH$, each intersecting the bag.
Note that induced $\{P_1\}$-packing treewidth is equivalent to tree-independence number,
while induced $\{P_2\}$-packing treewidth is equivalent to induced matching treewidth.
Let us look at this family of parameters through the lens of questions \ref{q1}, \ref{q2}, and \ref{q3} stated in the introduction. We focus on the first open case, i.e., $H = \{P_3\}$.

Question \ref{q1} was partially answered as boundedness of induced packing yields \emph{quasipolynomial-time} algorithms for some problems including \MWIS, \textsc{List 3-Coloring}, and \textsc{Odd Cycle Transversal}~\cite{amirwidth}.
Question \ref{q2} is answered as a simple argument allows us to reduce the problem to computing a tree decomposition of small independence number in an auxiliary graph, which yields an \textsf{XP}-approximation algorithm for induced $\{P_3\}$-packing treewidth.
However, \ref{q3} is wide open.
In particular, it would be interesting to investigate which classes (say, defined by a forbidding a single graph) have bounded induced $\{P_3\}$-packing treewidth.

\newpage

\bibliographystyle{plain}
\bibliography{references}

\end{document}